# UPPER BOUND ON THE DISCONNECTION TIME OF DISCRETE CYLINDERS AND RANDOM INTERLACEMENTS

By Alain-Sol Sznitman

*ETH Zürich*

We study the asymptotic behavior for large $N$ of the disconnection time $T_N$ of a simple random walk on the discrete cylinder $(\mathbb{Z}/N\mathbb{Z})^d \times \mathbb{Z}$, when $d \geq 2$. We explore its connection with the model of random interlacements on $\mathbb{Z}^{d+1}$ recently introduced in [*Ann. Math.*, in press], and specifically with the percolative properties of the vacant set left by random interlacements. As an application we show that in the large $N$ limit the tail of $T_N/N^{2d}$ is dominated by the tail of the first time when the supremum over the space variable of the Brownian local times reaches a certain critical value. As a by-product, we prove the tightness of the laws of $T_N/N^{2d}$, when $d \geq 2$.

**0. Introduction.** The present article derives an upper bound on the disconnection time $T_N$ of a discrete cylinder with base a $d$-dimensional torus of large side-length $N$. It explores some of the connections of this question with the percolative properties of the model of random interlacements recently introduced in [12]. A variety of results concerning the disconnection time by simple random walk of discrete cylinders with various large bases has been obtained; cf. [3, 4] and [11]. In particular, it appears that with broad generality the disconnection time has a rough order of magnitude comparable to the square of the number of points in the base. In all the above quoted references upper bounds on the disconnection time hinge on the fact that once the walk has covered the "zero level" of the cylinder, disconnection has occurred. This causes the appearance in the resulting upper bounds of spurious factors involving some power of the logarithm of the cardinality of the base. The present work departs from this approach and builds on the results of [12] concerning percolation for the vacant set left by random interlacements. Notably, we show here that when $d \geq 2$, the laws of $T_N/N^{2d}$ are tight. Together with the results of [4], this implies that when $d$ is sufficiently large, that is, $d \geq 17$, $T_N$ "lives in scale $N^{2d}$." Moreover, this work leads to a natural guess concerning the convergence and characterization of the limit of the distributions of $T_N/N^{2d}$ in terms of Brownian local times.

Before discussing these matters any further, let us first present the model more precisely. For $d \geq 2$ and $N \geq 1$ we consider the discrete cylinder

$$(0.1) \qquad E = \mathbb{T} \times \mathbb{Z} \qquad \text{where } \mathbb{T} = (\mathbb{Z}/N\mathbb{Z})^d.$$

A finite subset $S \subseteq E$ is said to disconnect $E$ if for large $M$, $\mathbb{T} \times [M, \infty)$ and $\mathbb{T} \times (-\infty, -M]$ are contained in distinct connected components of $E \setminus S$. For $x$ in $E$ we denote with $P_x$ the canonical law on $E^{\mathbb{N}}$ of a simple random walk on $E$ starting at $x$, and write $E_x$ for the corresponding expectation. We let $X. = (Y., Z.)$ stand for the canonical process with $Y.$ and $Z.$ its respective $\mathbb{T}$- and $\mathbb{Z}$-components. A key object of interest in this article is the disconnection time

$$(0.2) \qquad T_N = \inf\{n \geq 0; X_{[0,n]} \text{ disconnects } E\}.$$

We write $\rho_k$, $k \geq 0$, for the times of successive displacements of the "vertical" component $Z.$ of $X.$, that is, $\rho_0 = 0$ and $\rho_k = \inf\{n > \rho_{k-1}; Z_n \neq Z_{\rho_{k-1}}\}$, for $k \geq 1$, as well as $\widehat{Z}.$ for the time changed process

$$(0.3) \qquad \widehat{Z}_k = Z_{\rho_k}, \qquad k \geq 0,$$

which is distributed as a one-dimensional simple random walk. The local time of $\widehat{Z}$ is defined as

$$(0.4) \qquad \widehat{L}_k^z = \sum_{0 \leq m < k} 1\{\widehat{Z}_m = z\} \qquad \text{with } k \geq 0, z \in \mathbb{Z}.$$

We also consider $\gamma_v^z$, the first time when the number of distinct visits of the walk $X$ to $\mathbb{T} \times \{z\}$, the "level $z$" in the cylinder, reaches an amount $v$:

$$(0.5) \qquad \gamma_v^z = \inf\{\rho_k; k \geq 0 \text{ and } \widehat{L}_k^z \geq v\} \qquad \text{with } v \geq 0, z \in \mathbb{Z}.$$

A further ingredient is the so-called random interlacement at level $u \geq 0$, introduced in [12]. It is the trace left on $\mathbb{Z}^{d+1}$ (here $d+1$, with $d \geq 2$, plays the role of $d \geq 3$ in [12]) by a cloud of paths constituting a Poisson point process in the space of doubly infinite trajectories modulo time-shift, tending to infinity at positive and negative infinite times. We refer to Section 1 for precise definitions. The nonnegative parameter $u$ is in essence a multiplicative factor of the intensity measure of this point process. In a standard fashion one constructs on the same space $(\Omega, \mathcal{A}, \mathbb{P})$ [see (1.25) and below (1.26)] the family $\mathcal{I}^u$, $u \geq 0$, of random interlacements at level $u$; see (1.32). They are the traces on $\mathbb{Z}^{d+1}$ of the cloud of trajectories modulo time-shift, "up to level $u$." The random subsets $\mathcal{I}^u$ increase with $u$ and for $u > 0$, they are infinite random connected subsets of $\mathbb{Z}^{d+1}$, ergodic under space translations;

cf. Theorem 2.1, and Corollary 2.3 of [12]. The complement of $\mathcal{I}^u$, denoted by $\mathcal{V}^u$, is the so-called vacant set at level $u$; see (1.34). An important role is played here by the critical parameter

(0.6)
$$u_{**} = \inf\{u \geq 0, \alpha(u) > 0\} \quad \text{with}$$
$$\alpha(u) = \sup\Big\{\alpha \geq 0, \lim_{L \to \infty} L^\alpha \mathbb{P}[B(0,L) \xleftrightarrow{\mathcal{V}^u} S(0,2L)] = 0\Big\} \quad \text{for } u \geq 0,$$

where the supremum is by convention 0 if the set in the second line of (0.6) is empty, and $\{B(0,L) \xleftrightarrow{\mathcal{V}^u} S(0,2L)\}$ denotes the event where there is a nearest neighbor path in $\mathcal{V}^u$ starting in $B(0,L)$, the closed ball of radius $L$ and center 0 for the $\ell^\infty$-distance, and ending in $S(0,2L)$, the $\ell^\infty$-sphere with radius $2L$ and center 0. We show in Lemma 1.4 that, for all $d \geq 2$,

(0.7)
$$u_* \leq u_{**} < \infty,$$

where $u_*$ is the critical parameter introduced in [12], such that $\mathbb{P}$-a.s. $\mathcal{V}^u$ has an infinite connected component, that is, percolates, when $u < u_*$, and only finite components when $u > u_*$. Among the key results of [12] are the facts that when $d \geq 2$, $u_* < \infty$ (cf. Theorem 3.5 of [12]) and $u_* > 0$, at least when $d \geq 6$ (cf. Theorem 4.3 of [12]; we recall that here $d+1$ plays the role of $d$ in [12]). This has later been extended to all $d \geq 2$ in Theorem 3.4 of [10]. It is a natural question whether actually $u_* = u_{**}$.

The main results of this article relate $T_N$ to $u_{**}$. Specifically, we show in Theorem 4.1 that, when $d \geq 2$,

(0.8)
$$\lim_N P_0\Big[T_N > \inf_{z \in \mathbb{Z}} \gamma^z_{(N^d/(d+1))u}\Big] = 0 \quad \text{when } u > u_{**}.$$

Loosely speaking, this says that, given $u > u_{**}$, when $N$ is large, once the number of distinct visits of the walk to some level $\mathbb{T} \times \{z\}$ of the cylinder exceeds $\frac{N^d}{(d+1)}u$, then typically disconnection must have occurred. This result has some similar flavor to [13], where the trace left by random walk in the neighborhood of points of the cylinder by times of order $N^{2d}$ is compared to random interlacements. As a consequence of the key property (0.8), we show in Corollary 4.6 that

(0.9)
$$\overline{\lim_N} P_0[T_N \geq sN^{2d}] \leq W\Big[\zeta\Big(\frac{u_{**}}{\sqrt{d+1}}\Big) \geq s\Big] \quad \text{for all } s > 0,$$

where $W$ stands for the Wiener measure and

(0.10)
$$\zeta(u) = \inf\Big\{t \geq 0; \sup_{v \in \mathbb{R}} L(v,t) \geq u\Big\} \quad \text{for } u \geq 0,$$

with $L(v,t)$ a jointly continuous version of the local time of the canonical Brownian motion. In particular, this shows that the laws of $T_N/N^{2d}$ under

$P_0$ are tight. With the results of [4], it also proves that when $d$ is large, that is, $d \geq 17$,

(0.11)   the laws on $(0, \infty)$ of $T_N/N^{2d}$ under $P_0$, with $N \geq 2$, are tight,

that is, "$T_N$ lives in scale $N^{2d}$." It is natural to wonder whether more than (0.9) holds and $T_N/N^{2d}$ actually converges in distribution toward $\zeta(\frac{u_{**}}{\sqrt{d+1}}) \stackrel{\text{law}}{=} \frac{u_{**}^2}{d+1}\zeta(1)$; see Remark 4.7 (and a related question whether $u_* = u_{**}$). Let us also mention that, thanks to the works [1] and [5], the Laplace transform of $\zeta(u)$ is known and can be expressed as

$$(0.12) \quad E^W[e^{-\theta^2/2\zeta(u)}] = \frac{\theta u}{[\sinh(\theta u/2)]^2} \frac{I_1(\theta u/2)}{I_0(\theta u/2)} \quad \text{for } \theta, u > 0,$$

with $I_\nu$ the modified Bessel function of index $\nu$.

We will now briefly sketch the strategy of the proof of the main Theorem 4.1. The rough idea we use in order to show (0.8) is that once sufficiently many distinct visits of a given level $z$ of the cylinder have taken place, that is, more than $\frac{N^d}{d+1}(u_{**} + \delta)$ distinct visits, then the trace left by the walk in a box with center at level $z$ and side-length $N^{1-\varepsilon}$, where $\varepsilon$ can be chosen arbitrarily small, dominates the trace left by random interlacements at level $u_{**} + \delta'$ in such a box, where $0 < \delta' (= \frac{\delta}{8}) < \delta$. With a straightforward covering argument and the definition of the critical exponent $u_{**}$ [cf. (0.6)], one finds by adjusting parameters that the probability of existence of a nearest neighbor path in the cylinder between levels $z - N^{1-\varepsilon}$ and $z + N^{1-\varepsilon}$, avoiding the trajectory of the walk, tends to 0 as $N$ goes to infinity. In order to take care of the infimum over $z$ which appears in (0.8), the above rough scheme is combined with an argument relying on the spatial regularity of the local time of the simple random walk $\widehat{Z}$. It enables us to simply consider a large but finite number of levels, regularly spaced at heights which differ by a small multiple of $N^d$.

The above-mentioned scheme crucially involves a stochastic domination argument; see, in particular, Proposition 4.2 and its proof. Its implementation goes through several steps. It begins with the extraction of excursions of the walk, which roughly correspond to successive returns to the box $B(z) = \mathbb{T} \times [z - N, z + N]$ and departures from the box $\widetilde{B}(z) = \mathbb{T} \times (z - h_N, z + h_N)$, where the height $h_N = [N(\log N)^2]$ is large enough so that the $\mathbb{T}$-component of the walk has time to homogenize between various excursions of the walk. There is, however, a special recipe in the precise specification of the excursions [cf. (2.2)], and it plays an important role. With the coupling techniques developed in Proposition 2.2, we are able to replace the true excursions of the walk with a collection of i.i.d. excursions for which the starting point is uniformly distributed on the union of the two levels $\mathbb{T} \times \{z + N\}$ and $\mathbb{T} \times \{z - N\}$, and the path otherwise evolves as a simple random walk on $E$ stopped when exiting $\widetilde{B}(z)$. The specific choice of this starting distribution leads to a key identity for the entrance law of the excursion in a subset $A$ of $\mathbb{T} \times (z - N, z + N) \subseteq B(z)$; see Lemma 1.1.

Via a Poissonization argument, the above mentioned identity induces a very handy comparison of the trace left by a Poisson number of i.i.d. excursions in a box of $E$ with center in $\mathbb{T} \times \{z\}$ and size $N^{1-\varepsilon}$, and the trace left in a box of the same size by trajectories of a random interlacement at a suitably calibrated level $u$, when the trajectories entering the box are stopped once they leave a concentric box of side-length $\frac{N}{2}$; see (4.22). To handle the truncation involved in stopping trajectories, we use Theorem 3.1, which shows that in essence the trace of the truncated trajectories in the box of size $N^{1-\varepsilon}$ dominates the untruncated trace in the same box of a random interlacement at a slightly lower level $u'$ (which in the application in Section 4 can be chosen equal to $u_{**} + \frac{\delta}{8}$). The important Theorem 3.1 solely pertains to the model of random interlacements. Its proof uses a "sprinkling technique" with a similar flavor to some of the arguments employed in Section 3 of [12], when showing that $u_* < \infty$. These are some of the main ingredients entering the proof of Theorem 4.1.

Let us now describe how the article is organized.

In Section 1 we introduce further notation and recall various useful facts concerning random walks and random interlacements. The key identity of the entrance law in sets interior to $B(z)$ of the specially tailored excursions appears in Lemma 1.1. The finiteness of $u_{**}$ is shown in Lemma 1.4.

In Section 2 we develop the coupling technique, which enables us to work with i.i.d. excursions in the sequel. The main result appears in Proposition 2.2.

In Section 3 we develop the sprinkling technique which shows that the trace left in a box of size of order $N^{1-\varepsilon}$ by trajectories of an interlacement at level $u$ stopped when exiting a concentric box of side-length $\frac{N}{2}$, in essence dominates the trace left in the box of size of order $N^{1-\varepsilon}$ by an interlacement at a slightly lower level $u'$. The main result is Theorem 3.1.

In Section 4 we prove the key statement (0.8) in Theorem 4.1, and its consequence (0.9) in Corollary 4.6. The proof of Theorem 4.1 is split into Proposition 4.2, where the key domination argument shows that once at a given level $z$ in the cylinder, sufficiently many distinct visits have occurred, then with high probability disconnection of the cylinder has taken place, and into Proposition 4.3, where the spatial regularity of the local time of $\widehat{Z}$ is used to replace the infimum over all levels $z$, which appears in (0.8), with an infimum over a large but finite number of levels.

Finally the convention concerning constants we use in the text is the following. Throughout $c$ or $c'$ denote positive constants, which solely depend on $d$, with values changing from place to place. The numbered constants $c_0, c_1, \ldots$ are fixed and refer to the value at their first appearance in the text.

Dependance of constants on additional parameters appears in the notation, for instance, $c(\varepsilon)$ denotes a positive constant depending on $d$ and $\varepsilon$.

**1. Notation and some useful properties.** In this section we introduce additional notation and present some useful results concerning random walks and random interlacements. In particular, the key identity for the hitting distribution of the excursions of the walk on the cylinder appears in Lemma 1.1, and the proof of the finiteness of the critical value $u_{**}$ of (0.6) is presented in Lemma 1.4.

We write $\mathbb{N} = \{0, 1, 2, \ldots\}$ for the set of natural numbers. Given a nonnegative real number $a$, we write $[a]$ for the integer part of $a$. We let $|\cdot|$ and $|\cdot|_\infty$, respectively, stand for the Euclidean and $\ell^\infty$-distances on $\mathbb{Z}^{d+1}$ or for the corresponding distances induced on $E$. Throughout the article we assume $d \geq 2$. We say that two points of $\mathbb{Z}^{d+1}$ or $E$ are neighbors if their $|\cdot|$-distance equals 1. With $B(x, r)$ and $S(x, r)$ we denote the closed $|\cdot|_\infty$-ball and $|\cdot|_\infty$-sphere with radius $r \geq 0$ and center $x$ in $\mathbb{Z}^{d+1}$ or $E$. For $A, B$ subsets of $\mathbb{Z}^{d+1}$ or $E$ we write $A + B$ for the set of elements $x + y$ with $x$ in $A$ and $y$ in $B$, and $d(A, B) = \inf\{|x - y|_\infty; x \in A, y \in B\}$ for the mutual $\ell^\infty$-distance of $A$ and $B$; when $A = \{x\}$ is a singleton we write $d(x, B)$ for simplicity. We also write $U \subset\subset \mathbb{Z}^{d+1}$ or $U \subset\subset E$ to indicate that $U$ is a finite subset of $\mathbb{Z}^{d+1}$ or $E$. Given $U$ subset of $\mathbb{Z}^{d+1}$ or $E$, we denote with $|U|$ the cardinality of $U$, with $\partial U$ the boundary of $U$ and $\partial_{\text{int}} U$ the interior boundary of $U$:

$$
\begin{aligned}
\partial U &= \{x \in U^c; \exists x' \in U, |x - x'| = 1\}, \\
\partial_{\text{int}} U &= \{x \in U; \exists x' \in U^c, |x - x'| = 1\}.
\end{aligned}
\tag{1.1}
$$

The canonical shift on $E^{\mathbb{N}}$ is denoted with $(\theta_n)_{n \geq 0}$, that is, $\theta_n$ stands for the map from $E^{\mathbb{N}}$ into $E^{\mathbb{N}}$ such that $(\theta_n w)(\cdot) = w(\cdot + n)$ for $w \in E^{\mathbb{N}}$, and we write $(\mathcal{F}_n)_{n \geq 0}$ for the canonical filtration. Given a subset $U$ of $E$, we denote with $H_U, \widetilde{H}_U$ and $T_U$ the entrance time, the hitting time of $U$ and the exit time from $U$:

$$
\begin{aligned}
H_U &= \inf\{n \geq 0; X_n \in U\}, \qquad \widetilde{H}_U = \inf\{n \geq 1; X_n \in U\}, \\
T_U &= \inf\{n \geq 0; X_n \notin U\}.
\end{aligned}
\tag{1.2}
$$

In the case of a singleton $U = \{x\}$, we simply write $H_x$ or $\widetilde{H}_x$. We denote with $P_x^{\mathbb{Z}^{d+1}}$ the canonical law of a simple random walk on $\mathbb{Z}^{d+1}$ starting at $x$ and with $E_x^{\mathbb{Z}^{d+1}}$ the corresponding expectation. We otherwise keep the same notation as for the walk on $E$ concerning the canonical process, the canonical shift and other natural objects such as in (1.2).

Given $K \subset\subset \mathbb{Z}^{d+1}$ and $U \supseteq K$, the equilibrium measure and capacity of $K$ relative to $U$ are defined by

$$(1.3) \qquad e_{K,U}(x) = \begin{cases} P_x^{\mathbb{Z}^{d+1}}[\widetilde{H}_K > T_U], & \text{for } x \in K, \\ 0, & \text{for } x \notin K, \end{cases}$$

and

$$(1.4) \qquad \operatorname{cap}_U(K) = \sum_{x \in K} e_{K,U}(x) \qquad [\text{note that } \operatorname{cap}_U(K) \leq |K|].$$

The Green function of the walk killed outside $U$ is defined as

$$(1.5) \quad g_U(x, x') = E_x^{\mathbb{Z}^{d+1}}\left[\sum_{n \geq 0} 1\{X_n = x', n < T_U\}\right] \qquad \text{for } x, x' \in \mathbb{Z}^{d+1}.$$

When $U = \mathbb{Z}^{d+1}$, we drop $U$ from the notation in (1.3)–(1.5). The Green function is symmetric in its two variables and the probability to enter $K$ before exiting $U$ can be expressed as

$$(1.6) \quad P_x^{\mathbb{Z}^{d+1}}[H_K < T_U] = \sum_{x' \in \mathbb{Z}^{d+1}} g_U(x, x') e_{K,U}(x') \qquad \text{for } x \in \mathbb{Z}^{d+1}.$$

One also has the bounds

$$(1.7) \quad \sum_{x' \in K} g_U(x, x') \Big/ \sup_{y \in K} \sum_{x' \in K} g_U(y, x') \leq P_x^{\mathbb{Z}^{d+1}}[H_K < T_U]$$

$$\leq \sum_{x' \in K} g_U(x, x') \Big/ \inf_{y \in K} \sum_{x' \in K} g_U(y, x').$$

These inequalities, for instance, follow from the $L^1(P_x^{\mathbb{Z}^{d+1}})$-convergence of the bounded martingale $M_n = \sum_{x' \in K} g_U(X_{n \wedge H_K \wedge T_U}, x')$, $n \geq 0$, toward $1\{H_K < T_U\} \sum_{x' \in K} g_U(X_{H_K}, x')$.

In the case of the discrete cylinder $E$, when $U \subsetneq E$ is a strict subset of $E$, we define the corresponding objects just as in (1.3)–(1.5) with $P_x$ and $E_x$ in place of $P_x^{\mathbb{Z}^{d+1}}$ and $E_x^{\mathbb{Z}^{d+1}}$. We then have a similar identity and bounds as in (1.6) and (1.7).

We will sometimes find it useful to consider the continuous time random walks $\overline{X}_\cdot$, $\overline{Y}_\cdot$ and $\overline{Z}_\cdot$ on $E$, $\mathbb{T}$ and $\mathbb{Z}$ with respective jump rates equal to $2(d+1), 2d$ and $2$. We denote (with some abuse of notation) by $P_x$, $P_y^{\mathbb{T}}$ and $P_z^{\mathbb{Z}}$ the corresponding canonical laws starting at $x \in E$, $y \in \mathbb{T}$ and $z \in \mathbb{Z}$. We otherwise use notation such as $(\overline{\theta}_t)_{t \geq 0}$, $(\overline{\mathcal{F}}_t)_{t \geq 0}$ or $\overline{H}_U$ to refer to the natural continuous time objects. The continuous time walks are convenient because, on the one hand, the discrete skeleton of $\overline{X}_\cdot$ is distributed as the discrete time walk $X_\cdot$ and, on the other hand, for $x = (y, z) \in E$,

$$(1.8) \quad \text{under } P_y^{\mathbb{T}} \times P_z^{\mathbb{Z}}, (\overline{Y}_\cdot, \overline{Z}_\cdot) \text{ has the canonical law } P_x \text{ governing } \overline{X}_\cdot.$$

One should, however, note that the discrete time processes $Y$ and $Z$, respective $\mathbb{T}$- and $\mathbb{Z}$-projections of $X$ [cf. (0.2) above], are not distributed as the discrete skeletons of $\overline{Y}$ and $\overline{Z}$; indeed, they need not jump at each integer time.

As mentioned in the Introduction, we will consider certain concentric boxes in the cylinder $E$ and certain excursions of the walk related to these boxes. More precisely, we introduce the height scales

$$r_N = N < h_N = [N(2 + (\log N)^2)] \tag{1.9}$$

as well as the boxes in $E$ centered at level $z \in \mathbb{Z}$,

$$B(z) = \mathbb{T} \times (z + I) \subseteq \widetilde{B}(z) = \mathbb{T} \times (z + \widetilde{I}), \tag{1.10}$$

$$\text{where } I = [-r_N, r_N] \text{ and } \widetilde{I} = (-h_N, h_N).$$

When $z = 0$, we simply write $B$ and $\widetilde{B}$. We also introduce the probability $q$ which is the equidistribution on the union of levels $r_N$, and $-r_N$ in $E$:

$$q = \frac{1}{2N^d} \sum_{x \in \mathbb{T} \times \{-r_N, r_N\}} \delta_x. \tag{1.11}$$

We now come to an identity which will be applied to the entrance distribution in subsets of $B \setminus \partial_{\text{int}} B$ prior to exit from $\widetilde{B}$ for the walk in $E$ with starting distribution $q$. This identity plays a crucial role in comparing the trace left by the walk in the neighborhood of points of $B$ away from $\partial_{\text{int}} B$, with random interlacements. For the sake of clarity, we state the result in a slightly more general form than needed. We consider

$$\widetilde{a} > a > b > \widetilde{b} \quad \text{in } \mathbb{Z} \tag{1.12}$$

$$\text{with } \frac{a+b}{2} = \frac{\widetilde{a}+\widetilde{b}}{2}, 2h = \widetilde{a} - \widetilde{b}, 2r = a - b,$$

so $h + r$ and $h - r$ are integers, but $h$ and $r$ are possibly half-integers. We then define the probability

$$q_{a,b} = \frac{1}{2N^d} \sum_{x \in \mathbb{T} \times \{a,b\}} \delta_x. \tag{1.13}$$

For a measure $\mu$ on $E$ we write $P_\mu$ in place of $\sum_{x \in E} \mu(x) P_x$. We can now state the following.

LEMMA 1.1 ($N \geq 3$). *If* $U = \mathbb{T} \times (\widetilde{b}, \widetilde{a})$, *one then has*

$$\sum_{x' \in E} q_{a,b}(x') g_U(x', x) = (d+1) \frac{(h-r)}{N^d} \quad \text{for all } x \in \mathbb{T} \times [b, a]. \tag{1.14}$$

*Moreover, for $K \subseteq \mathbb{T} \times (b, a)$ one also has*

$$(1.15) \quad P_{q_{a,b}}[X_{H_K} = x, H_K < T_U] = (d+1)\frac{(h-r)}{N^d} e_{K,U}(x) \quad \text{for } x \in K.$$

PROOF. We begin with the proof of (1.14). With the help of the continuous time process $\overline{X}_\cdot$ and the symmetry of $g_U(\cdot, \cdot)$, we can write for $x = (y, z) \in \mathbb{T} \times [b, a]$

$$\sum_{x' \in E} q_{a,b}(x') g_U(x', x)$$

$$= \sum_{x' \in E} q_{a,b}(x') g_U(x, x')$$

$$= \sum_{x' \in E} q_{a,b}(x') 2(d+1) E_x\left[\int_0^{T_U} 1\{\overline{X}_t = x'\} dt\right]$$

$$\stackrel{(1.8),(1.13)}{=} \sum_{y' \in \mathbb{T}} \frac{(d+1)}{N^d} E_y^{\mathbb{T}}$$

$$(1.16)$$

$$\times E_z^{\mathbb{Z}}\left[\int_0^\infty 1\{\overline{Y}_t = y'\} 1\{\overline{Z}_t = a, \overline{T}_{(\widetilde{b},\widetilde{a})} > t\} dt \right.$$

$$\left. + \int_0^\infty 1\{\overline{Y}_t = y'\} 1\{\overline{Z}_t = b, \overline{T}_{(\widetilde{b},\widetilde{a})} > t\} dt\right]$$

$$= \frac{(d+1)}{N^d}\left(E_z^{\mathbb{Z}}\left[\int_0^{\overline{T}_{(\widetilde{b},\widetilde{a})}} 1\{\overline{Z}_t = a\} dt\right] + E_z^{\mathbb{Z}}\left[\int_0^{\overline{T}_{(\widetilde{b},\widetilde{a})}} 1\{\overline{Z}_t = b\} dt\right]\right)$$

$$= \frac{(d+1)}{2N^d}\left(\frac{P_z^{\mathbb{Z}}[H_a < T_{(\widetilde{b},\widetilde{a})}]}{P_a^{\mathbb{Z}}[\widetilde{H}_a > T_{(\widetilde{b},\widetilde{a})}]} + \frac{P_z^{\mathbb{Z}}[H_b < T_{(\widetilde{b},\widetilde{a})}]}{P_b^{\mathbb{Z}}[\widetilde{H}_b > T_{(\widetilde{b},\widetilde{a})}]}\right),$$

using again the link between the continuous time and discrete walk, as well as a classical identity for the Green function of the discrete walk in the last step. Using symmetry around $\frac{a+b}{2} = \frac{\widetilde{a}+\widetilde{b}}{2}$, we see that

$$(1.17) \quad P_a^{\mathbb{Z}}[\widetilde{H}_a > T_{(\widetilde{b},\widetilde{a})}] = P_b^{\mathbb{Z}}[\widetilde{H}_b > T_{(\widetilde{b},\widetilde{a})}] = \tfrac{1}{2}(h-r)^{-1} + \tfrac{1}{2}(h+r)^{-1}.$$

Moreover, we also have

$$(1.18) \quad P_z^{\mathbb{Z}}[H_a < T_{(\widetilde{b},\widetilde{a})}] + P_z^{\mathbb{Z}}[H_b < T_{(\widetilde{b},\widetilde{a})}] = \frac{z - \widetilde{b}}{h+r} + \frac{\widetilde{a} - z}{h+r} = \frac{2h}{h+r}.$$

Inserting these identities in the last line of (1.16), we find that

$$\sum_{x' \in E} q_{a,b}(x') g_U(x', x) = \frac{(d+1)}{2N^d}\left[\tfrac{1}{2}(h-r)^{-1} + \tfrac{1}{2}(h+r)^{-1}\right]^{-1} \frac{2h}{h+r}$$

$$= \frac{(d+1)}{N^d}(h-r).$$

This proves (1.14).

We now turn to the proof of (1.15). For $x \in K (\subseteq \mathbb{T} \times (b,a))$ we have

(1.19)
$$\begin{aligned}
P_{q_{a,b}}[X_{H_K} = x, H_K < T_U] \\
= \sum_{n \geq 1} P_{q_{a,b}}[X_n = x, n = T_{U \setminus K}] \\
= \sum_{n \geq 1} P_{q_{a,b}}[T_{U \setminus K} > n-1, X_1 \circ \theta_{n-1} = x] \\
\stackrel{\text{Markov}}{=} \sum_{x' \in E} \sum_{x'' : |x''-x|=1} q_{a,b}(x') g_{U \setminus K}(x', x'') \frac{1}{2(d+1)}.
\end{aligned}$$

Note that for $x', x''$ in $E$, the application of the strong Markov property at time $H_K$ in the formula corresponding to (1.5) yields that

(1.20) $\quad g_U(x'', x') = g_{U \setminus K}(x'', x') + E_{x''}[H_K < T_U, g_U(X_{H_K}, x')].$

Coming back to the last line of (1.19), using the symmetry of the Green functions as well as (1.14), we see that for $x'' \in \mathbb{T} \times [b, a]$

(1.21)
$$\begin{aligned}
\sum_{x' \in E} q_{a,b}(x') g_{U \setminus K}(x', x'') = (d+1) \frac{(h-r)}{N^d}(1 - P_{x''}[H_K < T_U]) \\
= (d+1) \frac{(h-r)}{N^d} P_{x''}[H_K > T_U].
\end{aligned}$$

Inserting this identity in the last line of (1.19), we find after the application of the Markov property at time 1 that

(1.22)
$$P_{q_{a,b}}[X_{H_K} = x, H_K < T_U] = (d+1) \frac{(h-r)}{N^d} P_x[\widetilde{H}_K > T_U]$$
$$\text{for } x \in K,$$

and this proves (1.15). $\square$

REMARK 1.2. In what follows the above lemma will be applied to the special case $\widetilde{a} = h_N$, $\widetilde{b} = -h_N$, $a = r_N$, $b = -r_N$ [see (1.9)], so that $q_{a,b} = q$ in (1.11), and $U = \mathbb{T} \times (\widetilde{b}, \widetilde{a}) = \widetilde{B}$; see below (1.10). If one considers the time of the last visit to $K \subseteq \widetilde{B}$ of the walk prior to the exit from $\widetilde{B}$,

$$L_K^{\widetilde{B}} = \sup\{n \geq 0, X_n \in K, n < T_{\widetilde{B}}\},$$

where the supremum is by convention equal to $-1$, when the set in parenthesis is empty, an application of the simple Markov property classically yields that, for $x \in K$,

$$
\begin{aligned}
P_q[X_{L_K^{\widetilde{B}}} = x, H_K < T_{\widetilde{B}}] &= \sum_{x' \in E} q(x') g_{\widetilde{B}}(x', x) P_x[\widetilde{H}_K > T_{\widetilde{B}}] \\
&= \sum_{x' \in E} q(x') g_{\widetilde{B}}(x', x) e_{K, \widetilde{B}}(x).
\end{aligned}
\tag{1.23}
$$

When $K \subseteq B \setminus \partial_{\mathrm{int}} B = \mathbb{T} \times (-r_N, r_N)$, Lemma 1.1 shows that this expression remains the same when $L_K^{\widetilde{B}}$ is replaced with $H_K$. In fact, for any nearest neighbor $\widetilde{B}$-valued path $\tau(n)$, $0 \le n \le N_\tau$, having its starting point $x_s = \tau(0)$ and its endpoint $x_e = \tau(N_\tau)$ in the support of $e_{K, \widetilde{B}}(\cdot)$, one has, with the help of (1.15) and the strong Markov property,

$$
\begin{aligned}
P_q[H_K < T_{\widetilde{B}}, (X_{H_K + \cdot})_{0 \le \cdot \le L_K^{\widetilde{B}} - H_K} = \tau] \\
= (d+1) \frac{(h_N - r_N)}{N^d} e_{K, \widetilde{B}}(x_s) \\
\times P_{x_s}[X_n = \tau(n), 0 \le n \le N_\tau] e_{K, \widetilde{B}}(x_e).
\end{aligned}
\tag{1.24}
$$

This identity has a strong flavor of (1.29) below and underlies the link between excursions with entrance distribution $q$ on the cylinder $E$ and random interlacements on $\mathbb{Z}^{d+1}$. This will play a crucial role in Section 4.

We now recall some notation and results from [12] concerning random interlacements. We denote with $W$ the space of doubly infinite nearest neighbor $\mathbb{Z}^{d+1}$-valued trajectories which tend to infinity at positive and negative infinite times, and with $W^*$ the space of equivalence classes of trajectories in $W$ modulo time-shift. The canonical projection from $W$ onto $W^*$ is denoted by $\pi^*$. We endow $W$ with its canonical $\sigma$-algebra $\mathcal{W}$ and denote by $X_n, n \in \mathbb{Z}$, the canonical coordinates.

We endow $W^*$ with $\mathcal{W}^* = \{A \subseteq W^*; (\pi^*)^{-1}(A) \in \mathcal{W}\}$, the largest $\sigma$-algebra on $W^*$ for which $\pi^* : (W, \mathcal{W}) \to (W^*, \mathcal{W}^*)$ is measurable. We also consider $W_+$ the space of nearest neighbor $\mathbb{Z}^{d+1}$-valued trajectories defined for nonnegative times and tending to infinity. We denote with $\mathcal{W}_+$ and $X_n, n \ge 0$, the canonical $\sigma$-algebra and the canonical process. Since $d \ge 2$, the simple random walk (on $\mathbb{Z}^{d+1}$) is transient and $W_+$ has full measure for any $P_x^{\mathbb{Z}^{d+1}}$, $x \in \mathbb{Z}^{d+1}$ [see above (1.3)], and we view whenever convenient the law of simple random walk on $\mathbb{Z}^{d+1}$ starting from $x$ as a probability on $(W_+, \mathcal{W}_+)$. We consider the space of point measures on $W^* \times \mathbb{R}_+$:

$$
\Omega = \Big\{ \omega = \sum_{i \ge 0} \delta_{(w_i^*, u_i)}, \text{ with } (w_i^*, u_i) \in W^* \times \mathbb{R}_+, i \ge 0, \text{ and }
$$

(1.25)
$$\omega(W_K^* \times [0,u]) < \infty, \text{ for any } K \subset\subset \mathbb{Z}^{d+1}, u \geq 0\Big\},$$

where for $K \subset\subset \mathbb{Z}^{d+1}$, $W_K^* \subseteq W^*$ is the subset of trajectories modulo time-shift which enter $K$,

(1.26)
$$W_K^* = \pi^*(W_K) \quad \text{and} \quad W_K = \{w \in W; \text{ for some } n \in \mathbb{Z}, X_n(w) \in K\}.$$

We endow $\Omega$ with the $\sigma$-algebra $\mathcal{A}$ generated by the evaluation maps $\omega \to \omega(D)$, where $D$ runs over the product $\sigma$-algebra $\mathcal{W}^* \times \mathcal{B}(\mathbb{R}_+)$. We denote with $\mathbb{P}$ the probability on $(\Omega, \mathcal{A})$ which is the Poisson point measure with intensity $\nu(dw^*)\, du$, giving finite mass to the sets $W_K^* \times [0, u]$, for $K \subset\subset \mathbb{Z}^{d+1}$, $u \geq 0$, where $\nu$ is the unique $\sigma$-finite measure on $(W^*, \mathcal{W}^*)$ such that, for any $K \subset\subset \mathbb{Z}^{d+1}$ (cf. Theorem 1.1 of [12]),

$$1_{W_K^*} \nu = \pi^* \circ Q_K \tag{1.27}$$

with $Q_K$ the finite measure on $W_K^0$, the subset of $W_K$ of trajectories which enter $K$ for the first time at time 0, such that, for $A, B$ in $\mathcal{W}_+$, $x \in \mathbb{Z}^{d+1}$

(1.28)
$$Q_K[(X_{-n})_{n\geq 0} \in A, X_0 = x, (X_n)_{n\geq 0} \in B]$$
$$= P_x^{\mathbb{Z}^{d+1}}[A|\widetilde{H}_K = \infty] e_K(x) P_x^{\mathbb{Z}^{d+1}}[B],$$

where we recall $e_K(\cdot)$ stands for the equilibrium measure of $K$; cf. (1.3) and below (1.5).

REMARK 1.3. It is also shown in Theorem 1.1 of [12] that, for $A, B \in \mathcal{W}_+$, $L_K(w)$ the time of the last visit of $K$ by the trajectory $w \in W_K^0$, and $\tau(n)$, $0 \leq n \leq N_\tau$, a finite nearest neighbor trajectory on $\mathbb{Z}^{d+1}$ with starting point $x_s = \tau(0)$ and endpoint $x_e = \tau(N_\tau)$, both in the support of $e_K(\subseteq \partial_{\text{int}} K)$,

(1.29)
$$Q_K[(X_{-n})_{n\geq 0} \in A, (X_\cdot)_{0\leq\cdot\leq L_K} = \tau, (X_{n+L_K})_{n\geq 0} \in B]$$
$$= P_{x_s}^{\mathbb{Z}^{d+1}}[A|\widetilde{H}_K = \infty] e_K(x_s) P_{x_s}^{\mathbb{Z}^{d+1}}[X_n = \tau_n, 0 \leq n \leq N_\tau]$$
$$\times e_K(x_e) P_{x_e}^{\mathbb{Z}^{d+1}}[B|\widetilde{H}_K = \infty].$$

In the case $A = B = \mathcal{W}_+$ the above formula has a very similar flavor to (1.24). It is also shown in Theorem 1.1 of [12] that $\nu$ is invariant under time reversal of trajectories in $W^*$ and under translation of trajectories by a constant vector.

Given $K \subset\subset \mathbb{Z}^{d+1}$, $u \geq 0$, one further defines on $(\Omega, \mathcal{A})$ the random point process with state space the set of finite point measures on $(W_+, \mathcal{W}_+)$:

$$\mu_{K,u}(\omega) = \sum_{i\geq 0} \delta_{(w_i^*)^{K,+}} 1\{w_i^* \in W_K^*, u_i \leq u\} \qquad \text{for } \omega = \sum_{i\geq 0} \delta_{(w_i^*, u_i)}, \tag{1.30}$$

where $(w^*)^{K,+}$ stands for the trajectory in $W_+$ which follows step by step $w^* \in W_K^*$ from the time it first enters $K$. One then has the fact that (cf. Proposition 1.3 of [12]), for $K \subset\subset \mathbb{Z}^{d+1}$, $u \geq 0$,

(1.31) $\mu_{K,u}$ is a Poisson point process on $(W_+, \mathcal{W}_+)$ with intensity measure $uP_{e_K}^{\mathbb{Z}^{d+1}}$,

where the notation is similar to below (1.13).

Given $\omega \in \Omega$, the interlacement at level $u \geq 0$ is the subset of $\mathbb{Z}^{d+1}$:

$$（1.32） \quad \mathcal{I}^u(\omega) = \bigcup_{u_i \leq u} \text{range}(w_i^*) \quad \text{if } \omega = \sum_{i \geq 0} \delta_{(w_i^*, u_i)},$$

where for $w^* \in W^*$, $\text{range}(w^*) = w(\mathbb{Z})$, for any $w \in W$ with $\pi^*(w) = w^*$. One readily sees that

$$（1.33） \quad \mathcal{I}^u(\omega) = \bigcup_{K \subset\subset \mathbb{Z}^{d+1}} \bigcup_{w \in \text{Supp} \, \mu_{K,u}(\omega)} w(\mathbb{N}).$$

The vacant set at level $u$ is then defined as

$$（1.34） \quad \mathcal{V}^u(\omega) = \mathbb{Z}^{d+1} \setminus \mathcal{I}^u(\omega) \quad \text{for } \omega \in \Omega, u \geq 0.$$

One has

$$（1.35） \quad \mathbb{P}[\mathcal{V}^u \supseteq K] = \exp\{-u \, \text{cap}(K)\} \quad \text{for all } K \subset\subset \mathbb{Z}^{d+1},$$

and this property leads to a characterization of the law $Q_u$ on $\{0,1\}^{\mathbb{Z}^{d+1}}$ of the random subset $\mathcal{V}^u$; see Remark 2.2 of [12]. As recalled in the Introduction, $Q_u$ is ergodic under spatial translations (cf. Theorem 2.1 of [12]), and for $u > 0$, $\mathcal{I}^u(\omega)$ is $\mathbb{P}$-a.s. an infinite connected subset of $\mathbb{Z}^{d+1}$; cf. Corollary 2.3 of [12]. To measure the percolative properties of $\mathcal{V}^u$, one introduces the nonincreasing function on $\mathbb{R}_+$,

(1.36) $\eta(u) = \mathbb{P}[0$ belongs to an infinite connected component of $\mathcal{V}^u]$,

and the critical value

$$（1.37） \quad u_* = \inf\{u \geq 0; \eta(u) = 0\} \in [0, \infty].$$

The main results of [12] show in Theorem 3.5 that $\mathcal{V}^u$ does not percolate for large $u$ and in Theorem 4.3 that, for $d \geq 6$ (recall we work here on $\mathbb{Z}^{d+1}$), $\mathcal{V}^u$ percolates for small $u$, that is,

$$（1.38） \quad u_* < \infty, \quad \text{and} \quad \text{for } d \geq 6, \quad u_* > 0.$$

We will now deduce from the controls derived in [12] the finiteness of the critical parameter $u_{**}$ introduced in (0.6).

Lemma 1.4.

$$u_* \leq u_{**} < \infty \tag{1.39}$$

*(it is a natural question whether $u_* = u_{**}$, see Remark 1.5 below).*

Proof. The left-hand side inequality is straightforward. Indeed, with similar notation as in (0.6), if $u > u_{**}$, then, for $L \geq 1$,

$$\eta(u) \leq \mathbb{P}[B(0,L) \xleftrightarrow{\mathcal{V}^u} S(0,2L)] \tag{1.40}$$

and the right-hand side tends to 0 with $L$ from the definition of $u_{**}$ in (0.6). Hence, $u \geq u_*$, and the left-hand side inequality of (1.39) follows by letting $u$ tend to $u_{**}$. We will now prove that $u_{**}$ is finite. Define for $L_0 > 1$ and $a = (100(d+1))^{-1}$ the sequence of length scales

$$L_{n+1} = \ell_n L_n \quad \text{where } \ell_n = 100[L_n^a], n \geq 0. \tag{1.41}$$

If we now introduce for $n \geq 0$

$$C^{(n)} = [0, L_n)^{d+1} \cap \mathbb{Z}^{d+1} \quad \text{and} \quad \widetilde{C}^{(n)} = \bigcup_i (iL_n + C^{(n)}), \tag{1.42}$$

where the union is over indexes $i$ in $\mathbb{Z}^{d+1}$ such that $d(C^{(n)}, iL_n + C^{(n)}) \leq 1$, in the notation from the beginning of this section, then with (3.16), (3.67) and (3.68) of [12], one can choose $L_0$ and $u > 0$ such that

$$\mathbb{P}[A^{u,n}] \leq cL_n^{-1} \quad \text{for all } n \geq 0,$$
$$\text{where } A^{u,n} = \{C^{(n)} \xleftrightarrow{\mathcal{V}^u} \partial_{\text{int}} \widetilde{C}^{(n)}\}. \tag{1.43}$$

When $L$ is large we can find a unique $n \geq 0$ such that $L_n \leq L < L_{n+1}$, and we can cover $B(0,L)$ by at most $c\ell_n^{d+1}$ possibly overlapping translates of $C^{(n)}$ contained in $B(0,L)$, with the corresponding translate of $\widetilde{C}^{(n)}$ included in $B(0, 2L)$. As a result of translation invariance of $Q_u$, we see that, for large $L$,

$$\begin{aligned}
\mathbb{P}[B(0,L) \xleftrightarrow{\mathcal{V}^u} B(0,2L)] &\leq c\ell_n^{(d+1)} \mathbb{P}[A^{u,n}] \leq c\ell_n^{(d+1)} L_n^{-1} \\
&\stackrel{(1.41)}{\leq} cL_n^{-(1-a(d+1))} \\
&\stackrel{(1.41)}{\leq} cL_{n+1}^{-(1-a(d+1))/(1+a)} \\
&\leq cL^{-(1-a(d+1))/(1+a)}.
\end{aligned} \tag{1.44}$$

This shows that in the notation of (0.6), $\alpha(u) \geq \frac{1-a(d+1)}{1+a} > 0$ and, hence, $u_{**} < \infty$. This completes the proof of (1.39). $\square$

REMARK 1.5. In the case of Bernoulli percolation it is well known that, in the sub-critical phase, the probability that the origin is connected by an open path to $S(0, L)$ decays exponentially with $L$; cf. Theorem 5.4, page 88 of [6]. So far no quantitative estimate for percolation in the vacant set of random interlacements showing, for instance, that $\alpha(u) > 0$ for $u > u_*$ is known. It is a natural question whether in fact $u_* = u_{**}$.

**2. The coupling construction.** In this section we introduce in (2.2) excursions of the walk in the cylinder $E$ which take place sometimes during the return to $B(z)$ and the departure from $\widetilde{B}(z)$; cf. (1.10). Due to translation invariance, we will only need to focus on the case $z = 0$ in the sequel. With the choice of $h_N$ in (1.9), the $\mathbb{T}$-component of the walk has enough time to homogenize between one excursion and the next. At the beginning of the next excursion the distribution of the location of the starting point of the path is close to $q$ in (1.11); cf. Lemma 2.1. This enables us to construct in Proposition 2.2 a coupling of the excursions of the path with a sequence of i.i.d. excursions with starting distribution $q$. Although simpler, this coupling has a similar flavor to what was needed in Section 3 of [13]. It will be very handy when comparing the percolative properties of the vacant set left by the walk on $E$ with that of the vacant set of random interlacements on $\mathbb{Z}^{d+1}$ in Section 4.

We begin with some notation. Given $z \in \mathbb{Z}$, we consider the stopping time $\sigma^z$ which is the first time when the walk visits one of the two levels $z \pm r_N$ after reaching level $z$:

$$(2.1) \qquad \sigma^z = H_{\mathbb{T} \times (z + \{-r_N, r_N\})} \circ \theta_{H_{\mathbb{T} \times \{z\}}} + H_{\mathbb{T} \times \{z\}}$$

as well as the successive times

$$(2.2) \qquad \begin{aligned} \sigma_0^z &= \sigma^z, \qquad \tau_0^z = T_{\widetilde{B}(z)} \circ \theta_{\sigma^z} + \sigma^z \quad \text{and} \quad \text{for } k \geq 0, \\ \sigma_{k+1}^z &= \sigma^z \circ \theta_{\tau_k^z} + \tau_k^z, \qquad \tau_{k+1}^z = T_{\widetilde{B}(z)} \circ \theta_{\sigma_{k+1}^z} + \sigma_{k+1}^z, \end{aligned}$$

so that $P_0$-a.s., $0 < \sigma_0^z < \tau_0^z < \cdots < \sigma_k^z < \tau_k^z < \cdots < \infty$.

When $z = 0$ we drop the superscript $z$ for simplicity. We begin with the following.

LEMMA 2.1 ($N \geq 1$). *For all $x' \notin \widetilde{B}$ and $x \in \mathbb{T} \times \{-r_N, r_N\}$ one has*

$$(2.3) \qquad |P_{x'}[X_\sigma = x] - q(x)| \leq cN^{-4d}.$$

PROOF. The argument has a similar flavor to what appears in the proof of Lemma 3.1 of [13]. With (1.8) and the fact that the discrete skeleton of $\overline{X}.$ is distributed as $X.$, we see that the distribution of $X_\sigma$ under $P_{x'}$, for $x' = (y', z')$, coincides with the distribution of $(\overline{Y}_{\overline{\sigma}}, \overline{Z}_{\overline{\sigma}})$ under $P_{y'}^{\mathbb{T}} \times P_{z'}^{\mathbb{Z}}$, if $\overline{\sigma}$ is the first time $\overline{Z}.$ reaches $\{-r_N, r_N\}$ after reaching 0. Thus, using reflection of the path after $H_{\mathbb{T}\times\{0\}}$, we see that, for $x = (y, r_N)$,

$$
\begin{aligned}
(2.4) \quad P_{x'}[X_\sigma = x] &= \tfrac{1}{2} P_{x'}[Y_\sigma = y] = \tfrac{1}{2} P_{y'}^{\mathbb{T}} \times P_{z'}^{\mathbb{Z}}[\overline{Y}_{\overline{\sigma}} = y] \\
&= \tfrac{1}{2} E_{z'}^{\mathbb{Z}}[\mu_{\overline{\sigma}}^{y'}(y)],
\end{aligned}
$$

where for $t \geq 0$ we have set $\mu_t^{y'}(\cdot) = P_{y'}^{\mathbb{T}}[\overline{Y}_t = \cdot]$.

Since $\overline{\sigma} \geq \overline{H}_0$, and $|z'| \geq h_N > N(\log N)^2$, standard estimates on the displacement of a simple random walk in continuous time on $\mathbb{Z}$ (see, for instance, (2.22) of [11]) show that

$$
(2.5) \quad P_{z'}^{\mathbb{Z}}[\overline{\sigma} \leq N^2(\log N)^2] \leq cN^{-4d}.
$$

We thus see coming back to (2.4) and (1.11) that

$$
\begin{aligned}
(2.6) \quad &|P_{x'}[X_\sigma = x] - q(x)| \\
&\leq \tfrac{1}{2} E_{z'}^{\mathbb{Z}}[|\mu_{\overline{\sigma}}^{y'}(y) - N^{-d}|, \overline{\sigma} > N^2(\log N)^2] + cN^{-4d}.
\end{aligned}
$$

Letting $\lambda_{\mathbb{T}}$ stand for the spectral gap of the walk $\overline{Y}.$ on $\mathbb{T}$ (cf. (1.8) of [11]), Lemma 1.1 in [11] states that, for $t \geq t_{\mathbb{T}} \stackrel{\text{def}}{=} \lambda_{\mathbb{T}}^{-1} \log(2|\mathbb{T}|)$, one has

$$
(2.7) \quad |\mu_t^{y'}(y) N^d - 1| \leq \tfrac{1}{2} \exp\{-(t - t_{\mathbb{T}})\lambda_{\mathbb{T}}\} \qquad \text{for } t \geq t_{\mathbb{T}}.
$$

One can see that $\lambda_{\mathbb{T}} \geq cN^{-2}$, for $N \geq 2$, and, hence, $t_{\mathbb{T}} \leq cN^2 \log(2N^d)$; see, for instance, the end of the proof of Lemma 3.1 of [13]. Coming back to (2.6), we see that, for $N \geq c$, the right-hand side is smaller than $cN^{-4d}$. The case of $x = (y, -r_N)$ is treated analogously and, adjusting constants, this completes the proof of Lemma 2.1. $\square$

We now come to the coupling construction which is the main object of this section.

PROPOSITION 2.2 ($N \geq 1$). *One can construct on an auxiliary probability space $(\widetilde{\Omega}, \widetilde{\mathcal{A}}, \widetilde{P})$ two sequences $X_\cdot^k, k \geq 1$, and $\widetilde{X}_\cdot^k, k \geq 1$, of $E$-valued processes such that*

(2.8) $X_\cdot^k, k \geq 1$, under $\widetilde{P}$ has same distribution as $X_{(\sigma_k + \cdot) \wedge \tau_k}, k \geq 1$, under $P_0$,

(2.9) $\widetilde{X}_\cdot^k, k \geq 1$, under $\widetilde{P}$ are independent and each distributed as $X_{\cdot \wedge T_{\widetilde{B}}}$, under $P_q$,

(2.10) $\widetilde{P}[X_\cdot^k \neq \widetilde{X}_\cdot^k] \leq cN^{-3d} \qquad \text{for } k \geq 1.$

PROOF. The distributions of $X_\sigma$ under $P_{x'}$, with $x' \in E$, and under $P_q$ are concentrated on $\mathbb{T} \times \{-r_N, r_N\}$. It follows from Lemma 2.1 that when $x'$ belongs to $\widetilde{B}^c$ their total variation distance is smaller than $cN^{d-4d} = cN^{-3d}$. With Theorem 5.2, page 19 of [8], we can construct for any $x'$ in $\widetilde{B}^c$ a probability $\rho_{x'}(dx, d\widetilde{x})$ on $E^2$ such that, under $\rho_{x'}$,

(2.11) the first component has the same distribution as $X_\sigma$ under $P_{x'}$,

(2.12) the second component has distribution $q$

and

(2.13) $$\rho_{x'}(\{x \neq \widetilde{x}\}) \leq cN^{-3d}.$$

Let us denote with $\mathcal{T}_E$ the countable set of $E$-valued trajectories which reach $\widetilde{B}^c$ after a finite time and are constant from then on, and are nearest neighbor prior to that time. The auxiliary space we consider is $\widetilde{\Omega} = (\mathcal{T}_E \times \mathcal{T}_E)^{[1,\infty)}$ endowed with the canonical $\sigma$-algebra $\widetilde{\mathcal{A}}$. We denote with $X_{\cdot}^k$, $\widetilde{X}_{\cdot}^k$, $k \geq 1$, the canonical coordinates on $\widetilde{\Omega}$. The probability $\widetilde{P}$ on $(\widetilde{\Omega}, \widetilde{\mathcal{A}})$ is constructed as follows. We introduce the kernel $R_{x'}$ from $\widetilde{B}^c$ to $\mathcal{T}_E \times \mathcal{T}_E$ such that, for $x'$ in $\widetilde{B}^c$, and $w, \widetilde{w}$ in $\mathcal{T}_E$,

(2.14) $$R_{x'}((w, \widetilde{w})) = \int_{\{x=\widetilde{x}\}} \rho_{x'}(dx, d\widetilde{x}) P_x[X_{\cdot \wedge T_{\widetilde{B}}} = w(\cdot) = \widetilde{w}(\cdot)]$$
$$+ \int_{\{x \neq \widetilde{x}\}} \rho_{x'}(dx, d\widetilde{x}) P_x[X_{\cdot \wedge T_{\widetilde{B}}} = w(\cdot)] P_{x'}[X_{\cdot \wedge T_{\widetilde{B}}} = \widetilde{w}(\cdot)].$$

In other words, under $R_{x'}$ the ordered pair of starting points of the two trajectories has distribution $\rho_{x'}$ and, conditionally on these starting points, when both points coincide the two trajectories coincide as well and evolve as the walk on $E$ stopped when exiting $\widetilde{B}$, and when the starting points differ the two trajectories evolve as independent copies of the walk stopped when exiting $\widetilde{B}$. We then construct $\widetilde{P}$ as the law of the Markov chain on $(\mathcal{T}_E \times \mathcal{T}_E)^{[1,\infty)}$ such that

(2.15) $(X_{\cdot}^1, \widetilde{X}_{\cdot}^1)$ has distribution $\sum_{x' \in \widetilde{B}^c} P_0[X_{\tau_0} = x'] R_{x'}$

and

(2.16) $R_{X_{T_{\widetilde{B}}}^k}$ is the conditional law of $(X_{\cdot}^{k+1}, \widetilde{X}_{\cdot}^{k+1})$ given $X_{\cdot}^{k'}, \widetilde{X}_{\cdot}^{k'}, 1 \leq k' \leq k$.

With (2.13) and (2.14), it is immediate that (2.10) holds. With (2.12) and (2.14) under any $R_{x'}$, $x' \in \widetilde{B}^c$, the second component is distributed as $X_{\cdot \wedge T_{\widetilde{B}}}$ under $P_q$, and (2.9) follows. On the other hand, under $R_{x'}$ the first component is distributed as $X_{(\sigma_0 + \cdot) \wedge \tau_0}$ under $P_{x'}$ and (2.8) is a consequence of the strong Markov property for the walk on $E$ and (2.2). $\square$

With the help of Proposition 2.2, we will be able to replace the excursions $X_{(\sigma_k+\cdot)\wedge \tau_k}$, $k \geq 1$, under $P_0$ by the collection of i.i.d. excursions $\widetilde{X}_\cdot^k$, $k \geq 1$, under $\widetilde{P}$ which have the same law as $X_{\cdot \wedge T_{\widetilde{B}}}$ under $\widetilde{P}$. Together with Lemma 1.1, this will facilitate the task of comparing the trace left by the excursions of the walk $X_\cdot$ in a sub-box $A$ of $B \setminus \partial_{\text{int}} B$ with center at level 0 and side-length of order $N^{1-\epsilon}$ with the trace left by a well calibrated random interlacement on $A$ (suitably identified to a subset of $\mathbb{Z}^{d+1}$).

**3. Truncation, sprinkling and random interlacements.** The object of this section is to develop a stochastic domination result showing that when $A$ and $\widetilde{C}$ are boxes in $\mathbb{Z}^{d+1}$ centered at the origin with respective side-length of order $N^{1-\epsilon}$ and $N$, then for large $N$ one can in essence dominate the trace on $A$ of the random interlacement at level $u'$ by the trace on $A$ left by all trajectories in the support of $\mu_{A,u}$ stopped at the exit time of $\widetilde{C}$, if $u$ is slightly bigger than $u'$. We refer to (1.30) for the notation. Thus, sprinkling, that is, choosing $u$ slightly bigger than $u'$, compensates the truncation of trajectories. Our main result Theorem 3.1 directly pertains to random interlacements and will play an important role in the next section when relating the critical parameter $u_{**}$ of (0.6) to the disconnection of the discrete cylinder by a simple random walk. We begin with some notation.

We consider $0 < \varepsilon < 1$ and denote with $A \subseteq \widetilde{C}$ the boxes in $\mathbb{Z}^{d+1}$:

$$(3.1) \qquad A = B\left(0, 2\left[\frac{N^{1-\varepsilon}}{8}\right]\right) \subseteq \widetilde{C} = B\left(0, \left[\frac{N}{4}\right]\right).$$

Given $u > 0$, we introduce for $\omega \in \Omega$ [cf. (1.25)] the truncated interlacement

$$(3.2) \qquad \mathcal{I}_{\widetilde{C}}^u(\omega) = \bigcup_{w \in \text{Supp}\, \mu_{A,u}(\omega)} w([0, T_{\widetilde{C}}]),$$

where the notation appears in (1.30). We will now compare when $u'$ is "sufficiently smaller" than $u$ the trace on $A$ of $\mathcal{I}^{u'}(\omega)$, the random interlacement at level $u'$ [cf. (1.32)] to the trace on $A$ of $\mathcal{I}_{\widetilde{C}}^u(\omega)$. Our main result is as follows.

THEOREM 3.1 ($d \geq 2$, $u > u' > 0$, $0 < \varepsilon < 1$). *For $N \geq c(\varepsilon)$, whenever*

$$(3.3) \qquad u \geq u' \exp\left\{\frac{c_0}{\varepsilon} e^{-\sqrt{\log N}}\right\},$$

*then there exist $\mathcal{I}^*$, $\overline{\mathcal{I}}$ random subsets of $A$ such that*

(3.4) $\qquad \mathcal{I}^{u'} \cap A = \mathcal{I}^* \cup \overline{\mathcal{I}},$

(3.5) $\qquad \mathcal{I}^*, \overline{\mathcal{I}}$ *are independent under* $\mathbb{P}$,

(3.6) $\qquad \mathbb{P}[\overline{\mathcal{I}} \neq \varnothing] \leq u' N^{-d},$

(3.7) $\qquad \mathcal{I}^*$ *is stochastically dominated by* $\mathcal{I}_{\widetilde{C}}^u \cap A$.

PROOF. We now define the integer $M$ and the subbox $C$ of $\widetilde{C}$ via

$$(3.8) \qquad M = [\exp\{\sqrt{\log N}\}] + 1, \qquad C = B\left(0, \left[\frac{N}{4M}\right]\right),$$

so that, for $N \geq c(\varepsilon)$,

$$(3.9) \qquad A \subseteq B(0, 100[N^{1-\varepsilon}]) \subseteq C \subseteq B\left(0, 100\left[\frac{N}{4M}\right]\right) \subseteq \widetilde{C}.$$

Throughout the proof we will write, for simplicity, $P_x$ and $E_x$ in place of $P_x^{\mathbb{Z}^{d+1}}$ and $E_x^{\mathbb{Z}^{d+1}}$, with $x$ in $\mathbb{Z}^{d+1}$, to denote the law on $(W_+, \mathcal{W}_+)$ of a simple random walk starting from $x$ and its corresponding expectation. We introduce the sequence of successive returns to $A$ and departures from $C$ of the walk, that is, with similar notation as in (1.2),

$$(3.10) \quad \begin{aligned} R_1 &= H_A, \qquad D_1 = T_C \circ \theta_{R_1} + R_1 \quad \text{and} \quad \text{for } k \geq 1, \\ R_{k+1} &= R_1 \circ \theta_{D_k} + D_k, \qquad D_{k+1} = D_1 \circ \theta_{D_k} + D_k, \end{aligned}$$

so that $0 \leq R_1 \leq D_1 \leq \cdots \leq R_k \leq D_k \leq \cdots \leq \infty$, and $P_x$-a.s. these inequalities, except maybe for the first one, are strict if the left-hand side is finite. Note that, for $\omega \in \Omega$ [see (1.25)], the finitely many trajectories of $W_+$ in the support of $\mu_{A,u'}(\omega)$ have a starting point in $\partial_{\text{int}} A \subseteq A$, and $D_1$ is finite for such trajectories. We can thus consider the index of the last finite exit from $C$ for the various trajectories in the support of $\mu_{A,u'}$ and write

$$(3.11) \quad \begin{aligned} \mu_{A,u'} &= \sum_{1 \leq \ell \leq r} \mu'_\ell + \overline{\mu} \qquad \text{where } r = \left[\frac{8}{\varepsilon}\right] + 1, \quad \text{and} \\ \mu'_\ell &= 1\{D_\ell < \infty = R_{\ell+1}\}\mu_{A,u'}, \qquad \overline{\mu} = 1\{D_{r+1} < \infty\}\mu_{A,u'}. \end{aligned}$$

Similarly, in the case of $\mu_{A,u}$ considering the last return to $A$ before exiting $\widetilde{C}$, we can write

$$(3.12) \qquad \mu_{A,u} = \sum_{\ell \geq 1} \mu_\ell \qquad \text{where } \mu_\ell = 1\{D_\ell < T_{\widetilde{C}} < R_{\ell+1}\}\mu_{A,u}.$$

As a direct consequence of (1.31) and the above decompositions, we see that under $\mathbb{P}$

$$(3.13) \quad \begin{aligned} &\mu'_\ell, 1 \leq \ell \leq r, \overline{\mu} \text{ are independent Poisson point processes on} \\ &(W_+, \mathcal{W}_+) \text{ with respective intensity measures } \zeta'_\ell = u'1\{D_\ell < \infty = \\ &R_{\ell+1}\}P_{e_A}, 1 \leq \ell \leq r, \text{ and } \overline{\zeta} = u'1\{D_{r+1} < \infty\}P_{e_A}, \end{aligned}$$

and that

$$(3.14) \quad \begin{aligned} &\mu_\ell, \ell \geq 1, \text{ are independent Poisson point processes on } (W_+, \mathcal{W}_+) \\ &\text{with respective intensity measures } \zeta_\ell = u1\{D_\ell < T_{\widetilde{C}} < R_{\ell+1}\}P_{e_A}. \end{aligned}$$

Moreover, we can express the respective traces of $\mathcal{I}^{u'}$ and $\mathcal{I}^{u}_{\widetilde{C}}$ on $A$ as follows:

(3.15)
$$\mathcal{I}^{u'} \cap A = \mathcal{I}^* \cup \overline{\mathcal{I}}, \qquad \text{where}$$
$$\mathcal{I}^* = \bigcup_{1 \leq \ell \leq r} \left( \bigcup_{w \in \operatorname{Supp} \mu'_\ell} w(\mathbb{N}) \cap A \right), \qquad \overline{\mathcal{I}} = \bigcup_{w \in \operatorname{Supp} \overline{\mu}} w(\mathbb{N}) \cap A$$

and

(3.16)
$$\mathcal{I}^{u}_{\widetilde{C}} \cap A = \bigcup_{\ell \geq 1} \left( \bigcup_{w \in \operatorname{Supp} \mu_\ell} w([0, T_{\widetilde{C}}]) \cap A \right).$$

Note that the successive application of the Markov property at times $D_r$, $D_{r-1}, \ldots, D_1$ yields for $N \geq c(\varepsilon)$

(3.17)
$$\overline{\zeta}(W_+) = u' P_{e_A}[R_{r+1} < \infty] \leq u' \left( \sup_{x \in \partial C} P_x[H_A < \infty] \right)^r \times \operatorname{cap}(A)$$
$$\leq u' \left\{ c \left( \frac{N^\varepsilon}{M} \right)^{-(d-1)} \right\}^r \times c N^{(1-\varepsilon)(d-1)}$$
$$\leq u' c^{r+1} N^{-3/4\varepsilon(d-1)r + (d-1)} \leq u' N^{-d},$$

where we have used the inequality in the right-hand side of (1.7) combined with standard bounds on the Green function (cf. [7], page 31) to estimate $\sup_{x \in \partial C} P_x[H_A < \infty]$, a standard upper bound on the capacity of $A$ (cf. (2.16), page 53 of [7]), the fact that $M$ grows slower than $N^{\varepsilon/4}$ [see (3.8)] and the definition of $r$ in (3.11).

We now introduce the measurable maps $\phi'_\ell$, for $\ell \geq 1$, from $\{D_\ell < \infty = R_{\ell+1}\}$ ($\subseteq W_+$) into $W_f^{\times \ell}$, where $W_f$ stands for the countable set of finite nearest neighbor trajectories on $\mathbb{Z}^{d+1}$ as well as the measurable maps $\phi_\ell$, $\ell \geq 1$, from $\{D_\ell < T_{\widetilde{C}} < R_{\ell+1}\}$ into $W_f^{\times \ell}$ defined through

(3.18)
$$\phi'_\ell(w) = (w(R_k + \cdot)_{0 \leq \cdot \leq D_k - R_k})_{1 \leq k \leq \ell} \qquad \text{for } w \in \{D_\ell < \infty = R_{\ell+1}\},$$
$$\phi_\ell(w) = (w(R_k + \cdot)_{0 \leq \cdot \leq D_k - R_k})_{1 \leq k \leq \ell} \qquad \text{for } w \in \{D_\ell < T_{\widetilde{C}} < R_{\ell+1}\}.$$

In other words, $\phi'_\ell(w)$, respectively, $\phi_\ell(w)$, keep track of the $\ell$ portions of the trajectory $w$ corresponding to times between the successive returns to $A$ up to the next departure from $C$. With (3.11) and (3.12), we can view $\mu'_\ell$ and $\mu_\ell$, for $\ell \geq 1$, as Poisson point processes on $\{D_\ell < \infty = R_{\ell+1}\}$ and $\{D_\ell < T_{\widetilde{C}} < R_{\ell+1}\}$, respectively. We denote with $\rho'_\ell$ and $\rho_\ell$ their respective images under the maps $\phi'_\ell$ and $\phi_\ell$. Hence, $\rho'_\ell$ and $\rho_\ell$ are Poisson point processes on $W_f^{\times \ell}$, and we write $\xi'_\ell$ and $\xi_\ell$ for their respective intensity. Note that, as a direct result of (3.13) and (3.14), we have

(3.19)  $\rho'_\ell, 1 \leq \ell \leq r$, and $\overline{\mu}$ are independent Poisson point processes,

(3.20)  $\rho_\ell, \ell \geq 1$, are independent Poisson point processes,

and, moreover, for $\ell \geq 1$,

$$\xi'_\ell(dw_1, \ldots, dw_\ell)$$
$$= u' P_{e_A}[D_\ell < R_{\ell+1} = \infty, (X_{R_k + \cdot})_{0 \leq \cdot \leq D_k - R_k} \in dw_k, 1 \leq k \leq \ell],$$
(3.21)
$$\xi_\ell(dw_1, \ldots, dw_\ell)$$
$$= u P_{e_A}[D_\ell < T_{\widetilde{C}} < R_{\ell+1}, (X_{R_k + \cdot})_{0 \leq \cdot \leq D_k - R_k} \in dw_k, 1 \leq k \leq \ell].$$

The next lemma will be useful in comparing $\xi'_\ell$ to $\xi_\ell$.

LEMMA 3.2 ($d \geq 2, 0 < \varepsilon < 1$). For $N \geq c(\varepsilon)$, one has for $x \in \partial C$ and $y \in \partial_{\mathrm{int}} A$

(3.22) $\quad P_x[T_{\widetilde{C}} < R_1 < \infty, X_{R_1} = y] \leq \dfrac{c_1}{M^{d-1}} P_x[R_1 < T_{\widetilde{C}}, X_{R_1} = y].$

PROOF. We implicitly assume (3.9). Note that for $y \in \partial_{\mathrm{int}} A$ one has

$$\sup_{z \in \partial C} P_z[T_{\widetilde{C}} < R_1 < \infty, X_{R_1} = y]$$

$$\leq \sup_{z \in \partial C} E_z[P_{X_{T_{\widetilde{C}}}}[R_1 < \infty, X_{R_1} = y]]$$
(3.23)
$$\leq \sup_{z \in \partial \widetilde{C}} P_z[H_{\partial C} < \infty] \sup_{z \in \partial C} P_z[R_1 < \infty, X_{R_1} = y]$$

$$\leq \dfrac{c}{M^{d-1}} \sup_{z \in \partial C} P_z[R_1 < \infty, X_{R_1} = y],$$

where in the last step we have used the rightmost inequality in (1.7) combined with standard bounds on the Green function just as in (3.17). Then observe that the function $z \to P_z[R_1 < \infty, X_{R_1} = y] = P_z[H_A < \infty, X_{H_A} = y]$ is positive harmonic on $A^c$. With the Harnack inequality (cf. [7], page 42) and a standard covering argument, we find that

(3.24) $\quad \sup_{z \in \partial C} P_z[R_1 < \infty, X_{R_1} = y] \leq c \inf_{z \in \partial C} P_z[R_1 < \infty, X_{R_1} = y].$

Therefore, coming back to (3.23), we see that

$$\sup_{z \in \partial C} P_z[T_{\widetilde{C}} < R_1 < \infty, X_{R_1} = y]$$

$$\leq \dfrac{c}{M^{d-1}} \inf_{\partial C} P_z[R_1 < \infty, X_{R_1} = y]$$
(3.25)
$$\leq \dfrac{c}{M^{d-1}} \inf_{\partial C} (P_z[T_{\widetilde{C}} < R_1 < \infty, X_{R_1} = y]$$
$$+ P_z[R_1 < T_{\widetilde{C}}, X_{R_1} = y]).$$

Assume that $N \geq c(\varepsilon)$ is such that $\frac{c}{M^{d-1}} \leq \frac{1}{2}$, with $c$ the constant appearing in the last member of (3.25), then one finds that, for $x \in \partial C$ and $y \in \partial_{\mathrm{int}} A$,

$$P_x[T_{\widetilde{C}} < R_1 < \infty, X_{R_1} = y] \leq \frac{2c}{M^{d-1}} P_x[R_1 < T_{\widetilde{C}}, X_{R_1} = y]$$

and this completes the proof of Lemma 3.2. $\square$

Our next step in the proof of Theorem 3.1 is the following.

LEMMA 3.3 ($d \geq 2, 0 \leq \varepsilon < 1$).    *For $N \geq c(\varepsilon)$, one has*

$$\xi'_\ell \leq \frac{u'}{u}\left(1 + \frac{c_1}{M^{d-1}}\right)^{\ell-1} \xi_\ell \qquad \text{for } \ell \geq 1. \tag{3.26}$$

PROOF.  With (3.21), we see that for $\ell \geq 1$, $w_1, \ldots, w_\ell \in W_f$, writing $w^s$ and $w^e$ for the respective starting point and endpoint of $w \in W_f$, one has

$$\xi'_\ell((w_1, \ldots, w_\ell))$$
$$= u' P_{e_A}[D_\ell < \infty = R_{\ell+1}, (X_{R_k + \cdot})_{0 \leq \cdot \leq D_k - R_k} = w_k(\cdot), 1 \leq k \leq \ell]$$
$$= \sum_{B \subseteq \{1,\ldots,\ell-1\}} u' P_{e_A}[D_\ell < T_{\widetilde{C}} < R_{\ell+1} = \infty,$$

(3.27)
$$(X_{R_k + \cdot})_{0 \leq \cdot \leq D_k - R_k} = w_k(\cdot), 1 \leq k \leq \ell,$$
$$T_{\widetilde{C}} \circ \theta_{D_k} + D_k < R_{k+1}, \text{ exactly when } k \in B,$$
$$\text{for } 1 \leq k \leq \ell - 1].$$

Note that the above expression vanishes unless $w^s_k \in \partial_{\mathrm{int}} A$, $w^e_k \in \partial C$, and $w_k$ takes values in $C$ except for its endpoint $w^e_k$, for each $1 \leq k \leq \ell$. If these conditions are fulfilled, we can use the strong Markov property repeatedly at times $D_\ell, R_\ell, D_{\ell-1}, \ldots, D_1$, and find that the last member of (3.27) equals

$$\sum_{B \subseteq \{1,\ldots,\ell-1\}} u' P_{e_A}[(X.)_{0 \leq \cdot \leq D_1} = w_1(\cdot)]$$
$$\times E_{w^e_1}[\mathbf{1}\{1 \notin B\}\mathbf{1}\{T_{\widetilde{C}} > R_1\} + \mathbf{1}\{1 \in B\}\mathbf{1}\{T_{\widetilde{C}} < R_1\},$$
$$R_1 < \infty, X_{R_1} = w^s_2]$$
$$\times P_{w^s_2}[(X.)_{0 \leq \cdot \leq D_1} = w_2(\cdot)] \cdots$$
$$\times E_{w^e_{\ell-1}}[\mathbf{1}\{\ell-1 \notin B\}\mathbf{1}\{T_{\widetilde{C}} > R_1\} + \mathbf{1}\{\ell-1 \in B\}\mathbf{1}\{T_{\widetilde{C}} < R_1\},$$
$$R_1 < \infty, X_{R_1} = w^s_\ell]$$

(3.28) $\quad \times P_{w^s_\ell}[(X.)_{0 \leq \cdot \leq D_1} = w_\ell(\cdot)] P_{w^e_\ell}[T_{\widetilde{C}} < R_1 = \infty]$

$$\stackrel{(3.22)}{\leq} \sum_{B \subseteq \{1,\ldots,\ell-1\}} \left(\frac{c_1}{M^{d-1}}\right)^{|B|} u' P_{e_A}[(X)_{0 \leq \cdot \leq D_1} = w_1(\cdot)]$$

$$\times P_{w_1^e}[R_1 < T_{\widetilde{C}}, X_{R_1} = w_2^s]$$
$$\times P_{w_2^s}[(X.)_{0 \leq \cdot \leq D_1} = w_2(\cdot)] \cdots$$
$$\times P_{w_{\ell-1}^e}[R_1 < T_{\widetilde{C}}, X_{R_1} = w_\ell^s]$$
$$\times P_{w_\ell^s}[(X.)_{0 \leq \cdot \leq D_1} = w_\ell(\cdot)] P_{w_\ell^e}[T_{\widetilde{C}} < R_1 = \infty].$$

Using the strong Markov property, we see that the above expression equals

$$u'\left(1 + \frac{c_1}{M^{d-1}}\right)^{\ell-1} P_{e_A}[T_{\widetilde{C}} \circ \theta_{D_k} + D_k > R_{k+1}, \text{ for } 1 \leq k \leq \ell-1,$$
$$(X_{R_k + \cdot})_{0 \leq \cdot \leq D_k - R_k} = w_k(\cdot), 1 \leq k \leq \ell,$$
$$D_\ell < T_{\widetilde{C}} \circ \theta_{D_\ell} + D_\ell < R_{\ell+1} = \infty]$$

(3.29)
$$\leq u'\left(1 + \frac{c_1}{M^{d-1}}\right)^{\ell-1} P_{e_A}[D_\ell < T_{\widetilde{C}} < R_{\ell+1},$$
$$(X_{R_k + \cdot})_{0 \leq \cdot \leq D_k - R_k} = w_k(\cdot), 1 \leq k \leq \ell]$$

$$\stackrel{(3.21)}{=} \frac{u'}{u}\left(1 + \frac{c_1}{M^{d-1}}\right)^{\ell-1} \xi_\ell((w_1,\ldots,w_\ell))$$

and this concludes the proof of Lemma 3.3. □

We now assume that

(3.30)
$$u \geq u' \exp\left\{\frac{8}{\varepsilon} \frac{c_1}{M^{d-1}}\right\}$$
$$\left[\geq u'\left(1 + \frac{c_1}{M^{d-1}}\right)^{\ell-1}, \text{ for all } 1 \leq \ell \leq r, \text{ see } (3.11)\right]$$

and find, as a consequence of Lemma 3.3, that

(3.31) $$\xi'_\ell \leq \xi_\ell \quad \text{for } 1 \leq \ell \leq r.$$

In view of (3.13) and (3.15), we see that

(3.32) $$\mathcal{I}^* \text{ and } \overline{\mathcal{I}} \text{ are independent under } \mathbb{P}$$

and that, with notation above (3.19),

(3.33) $$\mathcal{I}^* = \bigcup_{1 \leq \ell \leq r} \bigcup_{(w_1,\ldots,w_\ell) \in \text{Supp}\, \rho'_\ell} \text{range}\, w_1 \cup \cdots \cup \text{range}\, w_\ell.$$

We also see that, with (3.16),

$$(3.34) \quad \mathcal{I}_{\widetilde{C}}^u \cap A \supseteq \bigcup_{1 \leq \ell \leq r} \bigcup_{(w_1,\ldots,w_\ell) \in \mathrm{Supp}\,\rho_\ell} \mathrm{range}\, w_1 \cup \cdots \cup \mathrm{range}\, w_\ell.$$

In view of the independence stated in (3.19) and (3.20), and of the domination stated in (3.31), we see that under $\mathbb{P}$

$$(3.35) \qquad \mathcal{I}_{\widetilde{C}}^u \cap A \text{ stochastically dominates } \mathcal{I}^*.$$

Together with the fact that

$$(3.36) \qquad \mathbb{P}[\overline{\mathcal{I}} \neq 0] \leq \overline{\zeta}(W_+) \stackrel{(3.17)}{\leq} u' N^{-d}$$

and recalling (3.15) and (3.32), Theorem 3.1 now follows by choosing $c_0 = 8c_1$; see (3.30) and (3.8). $\square$

REMARK 3.4. It is clear from the proof of Theorem 3.1 that the specific choice of the factor $e^{-\sqrt{\log N}}$ inside the exponential in the right-hand side of (3.3) is not essential. One could just as well use a factor $1/\psi(N)$, where $\psi(\cdot)$ is a positive function on $[1,\infty)$ tending to infinity such that $\psi(t) = o(t^\gamma)$ for all $\gamma > 0$, and assuming $N \geq c(\varepsilon, \psi)$ in the statement of Theorem 3.1. The present choice will be sufficient for our purpose.

**4. Upper bound on the disconnection time.** We now come to the main object of the present article, namely, the derivation of the upper bound (0.8) on the disconnection time $T_N$ of the discrete cylinder $E$ (cf. Theorem 4.1) and its Corollary 4.6 relating the asymptotic behavior of $T_N$ to the Brownian stopping time $\zeta(\frac{u_{**}}{\sqrt{d+1}})$; see (0.9) and (0.10). The strategy employed to show Theorem 4.1 roughly goes as follows. We will show that once for some $z \in \mathbb{Z}$ the local time at $z$ of $\widehat{Z}$, [see (0.3)] exceeds $\frac{N^d}{(d+1)} u_0$ with $u_0 > u_{**}$, then typically all excursions $X_{[\sigma_k^z, \tau_k^z]}$, with $k \leq \frac{N^d}{(d+1)h_N} u_1$, have already occurred, where $u_{**} < u_1 < u_0$. In addition, an argument based on the spatial regularity of the local time will allow us to only consider a large but finite number of levels $z$'s in the cylinder as $N$ goes to infinity; see Proposition 4.3. With the coupling technique of Section 2, we will be able to replace the excursions $X_{[\sigma_k^z, \tau_k^z]}$, $1 \leq k \leq \frac{N^d}{(d+1)h_N} u_1$, by a collection of i.i.d. excursions with starting distribution, the vertical translation to level $z$ of $q$ in (1.11). With a Poissonization argument, it will suffice to consider a Poisson number of such i.i.d. excursions with parameter $\frac{N^d}{(d+1)h_N} u_2$, where $u_{**} < u_2 < u_1$. The special character of these excursions (see Lemma 1.1 and Remark 1.2) and the domination results for the trace of random interlacements of Section 3 will allow to compare the trace left by this Poisson number of excursions in a box of the cylinder with side-length $N^{1-\varepsilon}$ and center at level $z$, to the trace left by a random interlacement at level $u_3$, with $u_{**} < u_3 < u_2$, in a box of $\mathbb{Z}^{d+1}$ of the same side-length, where $\varepsilon$ will be chosen as a function of $\alpha(u_3)$, in the notation of (0.6). It will follow that disconnection of the cylinder typically must have occurred; see Proposition 4.2. Combining Propositions 4.2 and 4.3 will yield Theorem 4.1.

We recall the notation (0.5). Our main result is as follows.

THEOREM 4.1 ($d \geq 2$). *For any $\delta > 0$ one has*

(4.1) $$\lim_N P_0\Big[T_N > \inf_{z \in \mathbb{Z}} \gamma^z_{N^d/(d+1)(u_{**}+\delta)}\Big] = 0.$$

PROOF. We begin with a reduction step which shows that (4.1) is the consequence of two claims that will be subsequently proved in Propositions 4.2 and 4.3 below. Indeed, we can write for $L, N \geq 1$, in the notation of (2.2),

$$P_0\Big[T_N > \inf_{z \in \mathbb{Z}} \gamma^z_{N^d/(d+1)(u_{**}+\delta)}\Big]$$

(4.2) $$\leq P_0\Big[T_N > \inf_{z=\ell/LN^d, |\ell| \leq L^2} \tau^z_{[N^d/((d+1)h_N)(u_{**}+\delta/2)]}\Big]$$

$$+ P_0\Big[\inf_{z=\ell/LN^d, |\ell| \leq L^2} \tau^z_{[N^d/((d+1)h_N)(u_{**}+\delta/2)]} > \inf_{z \in \mathbb{Z}} \gamma^z_{N^d/(d+1)(u_{**}+\delta)}\Big].$$

As a result, we see that (4.1) will follow from the two propositions:

PROPOSITION 4.2 ($d \geq 2, \delta > 0$).

(4.3) *For all $z \in \mathbb{Z}$,* $\quad \lim_N P_0[T_N > \tau^z_{[N^d/((d+1)h_N)(u_{**}+\delta/2)]}] = 0.$

PROPOSITION 4.3 ($d \geq 2, \delta > 0$).

(4.4) $$\lim_L \overline{\lim_N} P_0\Big[\inf_{z=\ell/LN^d, |\ell| \leq L^2} \tau^z_{[N^d/((d+1)h_N)(u_{**}+\delta/2)]}$$
$$> \inf_{z \in \mathbb{Z}} \gamma^z_{N^d/(d+1)(u_{**}+\delta)}\Big] = 0.$$

We start with the following.

PROOF OF PROPOSITION 4.2. The application of the strong Markov property at the entrance time of the walk in $\mathbb{T} \times \{z\}$, together with translation invariance, shows that it suffices to consider the case $z = 0$ when proving (4.3). With (2.10) of Proposition 2.2, bringing the i.i.d. excursions $\widetilde{X}^k, k \geq 1$, into play, we see that (4.3) will follow once we show that

(4.5) $$\lim_N \widetilde{P}[\text{range}(\widetilde{X}^1_\cdot) \cup \cdots \cup \text{range}(\widetilde{X}^{[N^d/((d+1)h_N)(u_{**}+\delta/2)]}_\cdot)$$
$$\text{does not disconnect } E] = 0.$$

If we now introduce an independent Poisson random variable $K_\lambda$ with intensity

$$\lambda = \frac{N^d}{(d+1)h_N}\left(u_{**} + \frac{\delta}{4}\right), \tag{4.6}$$

then with a slight abuse of notation we have

$$\lim_N \widetilde{P}\left[K_\lambda > \left[\frac{N^d}{(d+1)h_N}\left(u_{**} + \frac{\delta}{2}\right)\right]\right] = 0$$

and, hence, the claim (4.3) follows from

(4.7) $\lim_N \widetilde{P}[\text{range}(\widetilde{X}^1_\cdot) \cup \cdots \cup \text{range}(\widetilde{X}^{K_\lambda}_\cdot) \text{ does not disconnect } E] = 0.$

We now choose [see (0.6) for the notation]

$$\varepsilon = \frac{1}{2d}\left(\alpha\left(u_{**} + \frac{\delta}{8}\right) \wedge 1\right) \in \left(0, \frac{1}{4}\right]. \tag{4.8}$$

We can cover $\mathbb{T} \times \{0\}$ by $cN^{\varepsilon d}$ closed $|\cdot|_\infty$-balls of radius $R = [\frac{N^{1-\varepsilon}}{8}]$ with center in $\mathbb{T} \times \{0\}$. Hence, using translation invariance, (4.7) follows from

(4.9) $\lim_N N^{\varepsilon d} \widetilde{P}[\text{there is a nearest neighbor path from } B(0,R) \text{ to } S(0,2R)$

$$\text{not intersecting range}(\widetilde{X}^1_\cdot) \cup \cdots \cup \text{range}(\widetilde{X}^{K_\lambda}_\cdot)] = 0.$$

We will write

$$A = B(0, 2R) \subseteq \widetilde{C} = B\left(0, \left[\frac{N}{4}\right]\right) \subseteq E \tag{4.10}$$

and for sufficiently large $N$, we will tacitly identify $\widetilde{C} \cup \partial \widetilde{C}$ with a subset of $\mathbb{Z}^{d+1}$, so that the notation agrees with (3.1). Given $\widetilde{X}^k_\cdot$, $k \geq 1$, entering $A$, we can define the nearest-neighbor trajectory $\overline{X}^k_\cdot$, which starts when $\widetilde{X}^k_\cdot$ enters $A$, follows $\widetilde{X}^k_\cdot$ and is stopped when $\widetilde{X}^k_\cdot$ exits $\widetilde{C}$. Then

$$\widetilde{\mu} = \sum_{1 \leq k \leq K_\lambda} 1\{\widetilde{X}^k_\cdot \text{ enters } A\} \delta_{\overline{X}^k_\cdot} \tag{4.11}$$

is a point process on the space of nearest neighbor $\widetilde{C} \cup \partial \widetilde{C}$-valued trajectories which are constant after a finite time. The key observation, in view of (1.15) of Lemma 1.1 when $U = \widetilde{B}$, (2.9) of Proposition 2.2 (and the main interest in introducing the independent Poisson variable $K_\lambda$), is that

$\widetilde{\mu}$ is a Poisson point process with intensity measure

(4.12) $\quad \frac{\lambda(d+1)}{N^d}(h_N - r_N)P_{e_{A,\widetilde{B}}}[X_{\cdot \wedge T_{\widetilde{C}}} \in \cdot]$

$$= \left(u_{**} + \frac{\delta}{4}\right)\left(1 - \frac{r_N}{h_N}\right)P_{e_{A,\widetilde{B}}}[X_{\cdot \wedge T_{\widetilde{C}}} \in \cdot].$$

We will now use the next lemma.

LEMMA 4.4 ($d \geq 2, \delta > 0$). *For $N \geq c(\delta)$, one has*

(4.13) $\quad$ *for all $x \in \partial_{\text{int}} A$,* $\quad e_{A,\widetilde{B}}(x) \geq e_A(x)\left(1 - c_2 \dfrac{(\log N)^2}{N^{(d-1)\varepsilon}}\right)$

*[see below (1.5) for the notation and recall $A \subseteq \widetilde{C}$ are viewed as subsets of both $E$ and $\mathbb{Z}^{d+1}$].*

PROOF. It is plain from (1.3) that, for $N \geq c$,

(4.14) $\qquad\qquad e_{A,\widetilde{C}}(x) \geq e_A(x) \qquad \text{for } x \in \partial_{\text{int}} A.$

It is therefore sufficient to prove (4.13) with $e_{A,\widetilde{C}}(x)$ in place of $e_A(x)$. On the other hand, with the analogue of (1.3) for the walk on $E$, we see that

(4.15)
$$\begin{aligned}
& e_{A,\widetilde{C}}(x) - e_{A,\widetilde{B}}(x) \\
&= P_x[T_{\widetilde{B}} > \widetilde{H}_A > T_{\widetilde{C}}] \\
&\overset{\text{strong Markov}}{=} E_x[\widetilde{H}_A > T_{\widetilde{C}}, P_{X_{T_{\widetilde{C}}}}[T_{\widetilde{B}} > H_A]] \\
&\leq e_{A,\widetilde{C}}(x) \sup_{x \in \partial\widetilde{C}} P_x[H_A < T_{\widetilde{B}}] \qquad \text{for } x \in \partial_{\text{int}} A.
\end{aligned}$$

Note that $\partial\widetilde{C} \subseteq S(0, [\tfrac{N}{4}]+1) \overset{\text{def}}{=} \overline{S}$, and the claim (4.13) will follow once we show that

(4.16) $\qquad \sup_{x \in \overline{S}} P_x[H_A < T_{\widetilde{B}}] \leq c \dfrac{(\log N)^2}{N^{(d-1)\varepsilon}} \qquad \text{for } N \geq c(\delta).$

To this end, consider the probability that the walk starting in $\overline{S}$ reaches $B(0, [\tfrac{1}{2} \times [\tfrac{N}{4}]])$ before hitting $\overline{S}$, and then enters $A$ before entering $\overline{S}$. We see with standard estimates on the one-dimensional simple random walk and the right-hand inequality of (1.7) combined with standard estimates on the Green function (cf. [7], page 31) that, for $N \geq c(\delta)$,

(4.17) $\quad \sup_{x \in \overline{S}} P_x[H_A < \widetilde{H}_{\overline{S}} \wedge T_{\widetilde{B}}] \leq cN^{-1}cN^{-(d-1)\varepsilon} = cN^{-1-(d-1)\varepsilon}.$

On the other hand, using estimates on the one-dimensional simple random walk to bound from below the probability to move at distance $[\tfrac{N}{10}]$ of $\widetilde{C} \cup \overline{S} = B(0,[\tfrac{N}{4}]+1)$ without hitting $\overline{S}$, the invariance principle to bound from below the probability to reach level $[\tfrac{N}{4}] + N$ in $E$ without entering $\overline{S}$, and once again estimates on the one-dimensional simple random walk to bound from below the probability to reach level $h_N$ before level $[\frac{N}{4}]+1$, we see that, for $N \geq c(\delta)$,

(4.18)
$$\inf_{x \in \overline{S}} P_x[T_{\widetilde{B}} < \widetilde{H}_{\overline{S}} \wedge H_A] \geq \frac{c}{N} \times c \times \frac{N-1}{h_N - [N/4] - 1}$$
$$\stackrel{(1.9)}{\geq} cN^{-1}(\log N)^{-2}.$$

We can then introduce the successive hitting times of $\overline{S}$,

(4.19) $$V_0 = 0, \qquad V_{k+1} = \widetilde{H}_{\overline{S}} \circ \theta_{V_k} + V_k, \qquad k \geq 0,$$

which are $P_x$-a.s. finite for all $x$ in $\overline{S}$ (and in $E$). Considering the pairwise disjoint events where $\theta_{V_m}^{-1}(\{H_A \wedge T_{\widetilde{B}} < \widetilde{H}_{\overline{S}}\})$, $m \geq 0$, first occurs when $m = k$, with $k \geq 0$, the application of the strong Markov property at time $V_k$ shows that, for $N \geq c(\delta)$, for all $x \in \overline{S}$,

$$P_x[H_A < T_{\widetilde{B}}]$$

(4.20)
$$\leq \frac{\sup_{x \in \overline{S}} P_x[H_A < \widetilde{H}_{\overline{S}} \wedge T_{\widetilde{B}}]}{\sup_{x \in \overline{S}} P_x[H_A < \widetilde{H}_{\overline{S}} \wedge T_{\widetilde{B}}] + \inf_{x \in \overline{S}} P_x[T_{\widetilde{B}} < \widetilde{H}_{\overline{S}} \wedge H_A]}$$
$$\stackrel{(4.17),\,(4.18)}{\leq} c(\log N)^2 N^{-(d-1)\varepsilon}.$$

This shows (4.16) and concludes the proof of Lemma 4.4. $\square$

We now proceed with the proof of (4.9). Note that with (4.11) one has

(4.21) $$(\mathrm{range}(\widetilde{X}_\cdot^1) \cup \cdots \cup \mathrm{range}(\widetilde{X}_\cdot^{K_\lambda})) \cap A \supseteq \bigcup_{w \in \mathrm{Supp}\,\widetilde{\mu}} (\mathrm{range}(w)) \cap A$$

and in view of (4.12) and (4.13), for $N \geq c(\delta)$,

(4.22) under $\widetilde{P}$, $(\mathrm{range}(\widetilde{X}_\cdot^1) \cup \cdots \cup \mathrm{range}(\widetilde{X}_\cdot^{K_\lambda})) \cap A$ stochastically dominates $\mathcal{I}_{\widetilde{C}}^u \cap A$ under $\mathbb{P}$ with

$$u = \left(u_{**} + \frac{\delta}{4}\right)\left(1 - \frac{r_N}{h_N}\right)\left(1 - c_2 \frac{(\log N)^2}{N^{(d-1)\varepsilon}}\right),$$

where we have used the fact stemming from (3.2) and (1.31) that

$$\mathcal{I}_{\widetilde{C}}^u(\omega) \cap A = \bigcup_{w \in \mathrm{Supp}\,\widetilde{\widetilde{\mu}}} (\mathrm{range}\,w) \cap A \qquad \text{where}$$

$$\widetilde{\widetilde{\mu}}(\omega) = \sum_{w \in \mathrm{Supp}\,\mu_{A,u}(\omega)} \delta_{w(\cdot \wedge T_{\widetilde{C}})} \text{ is a Poisson point process}$$

$$\text{with intensity measure } uP_{e_A}^{\mathbb{Z}^{d+1}}[X_{\cdot \wedge T_{\widetilde{C}}} \in \cdot].$$

Hence, returning to the expression in (4.9), we see that its lim sup over $N$ is smaller than [see (0.6) for notation]

$$\overline{\lim_{N}} N^{\varepsilon d}\mathbb{P}[\text{a nearest neighbor path in } (\mathcal{I}_C^u \cap A)^c$$

(4.23)
$$\text{joins } B(0,R) \text{ with } S(0,2R)].$$

If we now define

(4.24) $\quad u' = u \exp\left\{-\frac{c_0}{\varepsilon}e^{-\sqrt{\log N}}\right\} \geq u_{**} + \frac{\delta}{8} \quad \text{for } N \geq c(\delta),$

it follows from (3.7) that the above expression with a similar notation as in (0.6) is smaller than

$$\overline{\lim_{N}} N^{\varepsilon d}\mathbb{P}[B(0,R) \overset{(\mathcal{I}^*)^c}{\longleftrightarrow} S(0,2R)]$$

(4.25)
$$\overset{(3.4)}{\leq} \overline{\lim_{N}} N^{\varepsilon d}(\mathbb{P}[B(0,R) \overset{(\mathcal{V}^{u'})}{\longleftrightarrow} S(0,2R)] + \mathbb{P}[\overline{\mathcal{I}} \neq \phi])$$

$$\overset{(3.6),\,(4.24)}{\leq} \overline{\lim_{N}} N^{\varepsilon d}(\mathbb{P}[B(0,R) \overset{\mathcal{V}^{u_{**}+\delta/8}}{\longleftrightarrow} S(0,2R)] + u'N^{-d})$$

$$\overset{(0.6),\,(4.8)}{\leq} \overline{\lim_{N}} N^{\varepsilon d} N^{-\alpha(u_{**}+\delta/8)+\varepsilon} \overset{(4.8)}{=} 0.$$

This concludes the proof of (4.9), and hence of Proposition 4.2. □

Our next concern is Proposition 4.3.

PROOF OF PROPOSITION 4.3. Our first step is the following.

LEMMA 4.5 $(d \geq 2, \delta > 0)$.

(4.26)
$$\lim_{N} P_0[\tau^z_{[N^d/((d+1)h_N)(u_{**}+\delta/2)]} \geq \gamma^z_{N^d/((d+1))(u_{**}+3/4\delta)}] = 0$$

*for all $z \in \mathbb{Z}$.*

PROOF. Denote with $H_k^z$, $k \geq 1$, the successive times of entrance of $X$ at level $z$ after departure from $\widetilde{B}(z)$, that is,

(4.27)
$$H_0^z = H_{\mathbb{T} \times \{z\}} \quad \text{and}$$
$$H_{k+1}^z = H_{\mathbb{T} \times \{z\}} \circ \theta_{T_{\widetilde{B}(z)}} \circ \theta_{H_k^z} + T_{\widetilde{B}(z)} \circ \theta_{H_k^z} + H_k^z \quad \text{for } k \geq 0.$$

It follows from (2.2) that $\tau_k^z$ coincides with the exit time of $\widetilde{B}(z)$ after $H_k^z$:

(4.28) $\qquad \tau_k^z = T_{\widetilde{B}(z)} \circ \theta_{H_k^z} + H_k^z \quad \text{for } k \geq 0.$

Notice also that under $P_x$, for $x \in \mathbb{T} \times \{z\}$, the number of visits of $\widehat{Z}_\ell$, $\ell \geq 0$ [cf. (0.3)], to $z$ before exiting $z + \widetilde{I} = z + (-h_N, h_N)$, that is, $\sum_{\ell \geq 0} 1\{\widehat{Z}_\ell = z, \rho_\ell < T_{\widetilde{B}(z)}\}$ is distributed as a geometric variable with success probability $h_N^{-1}$. The application of the strong Markov property at the successive times $H_m^z$, $0 \leq m \leq k$, and (4.28) then shows that

(4.29) under $P_0$, $\sum_{\ell \geq 0} 1\{\widehat{Z}_\ell = z, \rho_\ell < \tau_k^z\}$ is distributed as the sum of $k+1$ independent geometric variables with success parameter $h_N^{-1}$.

Thus, choosing $k = [\frac{N^d}{(d+1)h_N}(u_{**} + \frac{\delta}{2})]$ and $\alpha = \frac{N^d}{(d+1)}(u_{**} + \frac{3}{4}\delta)$, we see that the probability which appears in (4.26) is equal to

$$P_0\left[\sum_{\ell \geq 0} 1\{\widehat{Z}_\ell = z, \rho_\ell < \tau_k^z\} \geq \alpha\right]$$

(4.30) $$\leq e^{-\lambda/h_N \alpha}\left(\frac{e^{\lambda/h_N}}{h_N} \frac{1}{1 - e^{\lambda/h_N}(1 - 1/h_N)}\right)^{k+1}$$

if $\lambda > 0$ and $e^{\lambda/h_N}\left(1 - \frac{1}{h_N}\right) < 1$,

where we have used (4.29) and the exponential Chebyshev inequality. If $\lambda < 1$ is small and fixed, for large $N$ the logarithm of the right member of (4.30) is equivalent to

$$-\frac{\lambda}{h_N}\frac{N^d}{(d+1)}\left(u_{**} + \frac{3}{4}\delta\right) + \frac{N^d}{(d+1)h_N}\left(u_{**} + \frac{\delta}{2}\right)\log\left(\frac{1}{1-\lambda}\right),$$

and this expression tends to $-\infty$. This concludes the proof of (4.26). $\square$

With Lemma 4.5, we see that, for given $L \geq 1$, the limsup over $N$ of the probability in (4.4) is bounded above by the lim sup over $N$ of the corresponding probability where $\tau_{[N^d/((d+1)h_N)(u_{**}+\delta/2)]}^z$ is replaced by $\gamma_{[N^d/(d+1)(u_{**}+3/4\delta)]}^z$. Hence, the claim (4.4) will follow once we show that

(4.31) $$\lim_L \overline{\lim}_N P_0\left[\inf_{z=\ell/LN^d, |\ell| \leq L^2} \gamma_{N^d/(d+1)(u_{**}+3/4\delta)}^z > \inf_{z \in \mathbb{Z}} \gamma_{N^d/(d+1)(u_{**}+\delta)}^z\right] = 0.$$

If we now introduce an integer $K \geq 1$, and write $\widetilde{\lim}$ in place of $\overline{\lim}_K \overline{\lim}_L \overline{\lim}_N$, we see that the above expression is smaller than

$$\widetilde{\lim} P_0\left[\rho_{KN^{2d}} \geq \inf_{z=\ell/LN^d, |\ell| \leq L^2} \gamma_{N^d/(d+1)(u_{**}+3/4\delta)}^z > \inf_{z \in \mathbb{Z}} \gamma_{N^d/(d+1)(u_{**}+\delta)}^z\right]$$

$$+ \widetilde{\lim} P_0[\gamma_{N^d/(d+1)(u_{**}+3/4\delta)}^0 > \rho_{KN^{2d}}]$$

(4.32)
$$\leq \widetilde{\lim} P_0 \left[ \sup_{k \leq KN^{2d}} \sup_{z \in \mathbb{Z}} \inf_{z' \in \{\ell/LN^d; |\ell| \leq L^2\}} |\widehat{L}_k^z - \widehat{L}_k^{z'}| \geq \frac{N^d}{(d+1)} \frac{\delta}{8} \right]$$
$$+ \widetilde{\lim_K} \overline{\lim_N} P_0 \left[ \widehat{L}_{KN^{2d}}^0 < \frac{N^d}{(d+1)} (u_{**} + 3/4\delta) \right].$$

With (1.20) of [2], one can construct on an auxiliary probability space $(\overline{\Omega}, \overline{\mathcal{A}}, \overline{P})$ a coupling of the local time $\widehat{L}$ of $\widehat{Z}$ with a jointly continuous version $L(\cdot, \cdot)$ of the Brownian local time so that

(4.33) $\overline{P}$-a.s., $\sup_{z \in \mathbb{Z}, k \geq 1} \frac{|\widehat{L}_k^z - L(z,k)|}{k^{1/4+\eta}} < \infty$ for all $\eta > 0$.

As a result, we see that the last member of (4.32) is smaller than

$$\widetilde{\lim} \overline{P} \left[ \sup_{t \leq KN^{2d}} \sup_{z \in \mathbb{Z}} \inf_{z' \in \{\ell/LN^d; |\ell| \leq L^2\}} |L(z,t) - L(z',t)| \geq \frac{N^d}{(d+1)} \frac{\delta}{16} \right]$$
$$+ \overline{\lim_K} \overline{\lim_N} \overline{P} \left[ L(0, KN^{2d}) \leq \frac{N^d}{(d+1)} (u_{**} + \delta) \right]$$
$$\stackrel{\text{scaling}}{\leq} \overline{\lim_K} \overline{\lim_L} \overline{P} \left[ \sup_{s \leq K} \sup_{v \in \mathbb{R}} \inf_{v' \in \{\ell/L, |\ell| \leq L^2\}} |L(v,s) - L(v',s)| \geq \frac{\delta}{16(d+1)} \right]$$
$$+ \overline{\lim_K} \overline{P} \left[ L(0, K) \leq \frac{u_{**} + \delta}{d+1} \right].$$

Since $\lim_{s \to \infty} L(0, s) = \infty$, $\overline{P}$-a.s., the last term vanishes, and since $\overline{P}$-a.s. the restriction to $\mathbb{R} \times [0, K]$ of $L(v, s)$ is continuous and compactly supported, the lim sup over $L$ of the probability in the previous line equals 0. Combining our estimates, we see that we have shown (4.31), and hence Proposition 4.3. □

As mentioned above (4.3), with Propositions 4.2 and 4.3, coming back to (4.2), we see that we have proved (4.1). This completes the proof of Theorem 4.1. □

As an application of Theorem 4.1, we will now derive an upper bound on $T_N$, which will, in particular, show that the variables $T_N/N^{2d}$ are tight. We recall from (0.10) the notation

$$\zeta(u) = \inf \left\{ t \geq 0; \sup_{v \in \mathbb{R}} L(v, t) \geq u \right\} \quad \text{for } u \geq 0,$$

with $L(\cdot, \cdot)$ a jointly continuous version of the local time of the canonical Brownian motion. Denoting with $W$ the Wiener measure, one has the scaling

property:

(4.34)    for $u \geq 0$,    $\zeta(u)$ and $u^2\zeta(1)$ have same law under $W$.

With [1] or [5] (cf. Proposition 5, page 89), as recalled in (0.12), one knows that, for $\theta, u \geq 0$,

$$(4.35) \qquad E^W[e^{-\theta^2/2\zeta(u)}] = \frac{\theta u}{(\sinh(\theta u/2))^2} \frac{I_1(\theta u/2)}{I_0(\theta u/2)}$$

with $I_\nu$ the modified Bessel function of order $\nu$; cf. [9], page 60.

COROLLARY 4.6 ($d \geq 2$).

$$(4.36) \quad \text{For } \gamma > 0, \qquad \varlimsup_N P_0[T_N \geq \gamma N^{2d}] \leq W\left[\zeta\left(\frac{u_{**}}{\sqrt{d+1}}\right) \geq \gamma\right],$$

and, in particular, the laws of $T_N/N^{2d}$ are tight.

PROOF.    Consider $0 < \gamma' < \gamma$, and $\delta > 0$. With Theorem 4.1, we see that

$$(4.37) \quad \varlimsup_N P_0[T_N \geq \gamma N^{2d}] \leq \varlimsup_N P_0\left[\inf_{z \in \mathbb{Z}} \gamma^z_{N^d/(d+1)(u_{**}+\delta)} \geq \gamma N^{2d}\right].$$

When $N \geq 3$, the sequence $\rho_k$, $k \geq 0$ [cf. below (0.2)], has the same distribution under $P_0$ as the partial sums of independent geometric variables with success probability $\frac{1}{d+1}$ (this distribution is independent of $N$). It follows from the strong law of large numbers that $P_0$-a.s., $\lim_k \frac{\rho_k}{k} = d+1$ and, hence, the right-hand side of (4.37) is smaller than

$$\varlimsup_N P_0\left[\inf_{z \in \mathbb{Z}} \gamma^z_{N^d/(d+1)(u_{**}+\delta)} > \rho_{[\gamma'/(d+1)N^{2d}]}\right]$$

$$\leq \varlimsup_N P_0\left[\sup_{z \in \mathbb{Z}} \widehat{L}^z_{[\gamma'/(d+1)N^{2d}]} < \frac{N^d}{d+1}(u_{**}+\delta)\right]$$

$$(4.38) \qquad \stackrel{(4.33)}{\leq} \varlimsup_N W\left[\sup_{z \in \mathbb{Z}} L\left(z, \left[\frac{\gamma'}{d+1}N^{2d}\right]\right) < \frac{N^d}{d+1}(u_{**}+2\delta)\right]$$

$$\stackrel{\text{scaling}}{=} \varlimsup_N W\left[\sup_{z \in \mathbb{Z}} L\left(\frac{z}{N^d}, \left[\frac{\gamma'}{d+1}N^{2d}\right]/N^{2d}\right) < \frac{u_{**}+2\delta}{d+1}\right]$$

$$\stackrel{\text{continuity}}{\leq} W\left[\sup_{v \in \mathbb{R}} L\left(v, \frac{\gamma'}{d+1}\right) \leq \frac{u_{**}+2\delta}{d+1}\right].$$

Letting $\gamma'$ tend to $\gamma$ and $\delta$ tend to 0, the above expression tends to

$$W\left[\sup_{v \in \mathbb{R}} L\left(v, \frac{\gamma}{d+1}\right) \leq \frac{u_{**}}{d+1}\right] \stackrel{\text{scaling}}{=} W\left[\sup_{v \in \mathbb{R}} L(v, \gamma) \leq \frac{u_{**}}{\sqrt{d+1}}\right].$$

One also knows (cf. [5], page 89 above Proposition 5) that, for $u \geq 0$,

(4.39) $\qquad W\text{-a.s.}, \qquad \zeta(u) = \inf\Big\{t \geq 0; \sup_{v \in \mathbb{R}} L(v,t) > u\Big\}$

and, therefore, the above expression equals $W[\zeta(\frac{u_{**}}{\sqrt{d+1}}) \geq \gamma]$, and this is an upper bound on the left-hand side of (4.37). This concludes the proof of Corollary 4.6. $\square$

REMARK 4.7. (1) Combined with the results of [4], Corollary 4.6 implies that when $d$ is large enough, that is, $d \geq 17$, the laws of $T_N/N^{2d}$ under $P_0$ with $N \geq 2$ are tight on $(0, \infty)$; see also (0.11).

(2) A natural question stemming from the present work is whether in fact

(4.40) $\qquad T_N/N^{2d}$ converges in distribution to $\zeta(\frac{u_*}{\sqrt{d+1}})$ as $N \to \infty$.

This question should be complemented by the further question whether it also holds that

(4.41) $$u_{**} = u_*$$

(one knows that $0 < u_* < \infty$ for $d+1 \geq 3$ (cf. [10] and [12]), and that $u_* \leq u_{**} < \infty$, for $d+1 \geq 3$, as shown in Lemma 1.4). These are just a few examples of the natural questions pertaining to the interplay between disconnection by a random walk of discrete cylinders and percolation for the vacant set of random interlacements.

DEPARTEMENT MATHEMATIK
ETH ZÜRICH
CH-8092 ZÜRICH
SWITZERLAND